\documentclass{article}
 
\usepackage{amsmath,amssymb,amsfonts,amsthm}
\usepackage{pstcol}

\def\ldot{\cdot_{\ltimes}}

\def\<{\langle}
\def\>{\rangle}

\def\ll{\lambda}
\def\dd{{\ \mathrm{d}}}
\def\Id{\mathrm{Id}}
\def\id{\mathrm{id}}
\def\Hom{\mathrm{Hom}}
\def\res{\langle z^{-1} \rangle}
\newcommand{\resz}[1]{\langle {#1}^{-1} \rangle}
\def\ix#1{_{(#1)}}

\def\C{\mathbb{C}}
\def\N{\mathbb{N}}
 
\def\A{{\cal{A}}}
\def\M{{\cal{M}}}

\def\Gi{G^{\mathrm{inv}}}
\def\Gd{G^{\mathrm{dif}}}

\def\Hd{{\cal{H}}^{\mathrm{dif}}}
\def\Hdt{{\cal{H}}^{\mathrm{dif}}_\ast}
\def\Hi{{\cal{H}}^{\mathrm{inv}}}
\def\Hit{{\cal{H}}^{\mathrm{inv}}_\ast}
\def\Hp{{\cal{H}}^{\gamma}}
\def\He{{\cal{H}}^{e}}
\def\Ha{{\cal{H}}^{\alpha}}
\def\Hat{\widetilde{\Ha}}

\def\B{{\cal{B}}}
\def\Bd{{\cal{B}}^{\mathrm{dif}}}
\def\BFdB{{\cal{B}}^{\mathrm{FdB}}}
\def\CC{{\cal{C}}}

\def\D{\Delta}
\def\Dd{\Delta^{\mathrm{dif}}}
\def\Ddt{\Delta^{\mathrm{dif}}_\ast}
\def\Di{\Delta^{\mathrm{inv}}}
\def\Dit{\Delta^{\mathrm{inv}}_\ast}
\def\Dl{\Delta^{\ltimes}}
\def\DB{\Delta^{\B}}
\def\DPp{\Delta^p_{\gamma}}
\def\DPe{\Delta^p_e}
\def\Da{\Delta^{\alpha}}

\def\Sd{S}

\def\ddif{\delta^{\mathrm{dif}}}
\def\da{\delta^{\alpha}}
\def\ddp{\delta^{\gamma}}
\def\dde{\delta^e}

\def\over{\slash}
\def\under{\backslash}

\newcommand{\lgraft}[2]
{\setlength{\unitlength}{1pt}
\begin{picture}(15,17)
        \put(11,0){#1}
        \put(2,15){#2}
        \put(6,8){$\backslash$}
        \end{picture}}
\newcommand{\rgraft}[2]
{\setlength{\unitlength}{1pt}
\begin{picture}(15,17)
        \put(2,0){#1}
        \put(11,15){#2}
        \put(6,8){$\slash$}
        \end{picture}}

\newtheorem{theorem}{Theorem}[section]
\newtheorem{proposition}[theorem]{Proposition}
\newtheorem{corollary}[theorem]{Corollary}
\newtheorem{lemma}[theorem]{Lemma}
\newtheorem{defin}[theorem]{Definition}
\newenvironment{definition}{\begin{defin} \em}{\end{defin}}
\newtheorem{rem}[theorem]{Remark}
\newenvironment{remark}{\begin{rem} \em}{\end{rem}}
\newenvironment{proof of}[1]{\bigskip \noindent{\em Proof of (\ref{#1}).\/}}
        {\hfill$\square$\par\vspace{.2cm}}
\numberwithin{equation}{section}
\numberwithin{figure}{section}

\renewcommand{\|}{
\setlength{\unitlength}{3pt}
\psset{unit=3pt}
\psset{runit=2pt}
\psset{linewidth=0.2}
\begin{pspicture}(0,.5)(2.5,4)
\psline(1,0)(1,3)
\end{pspicture}}

\newcommand{\Y}{
\setlength{\unitlength}{3pt}
\psset{unit=3pt}
\psset{runit=2pt}
\psset{linewidth=0.2}
\begin{pspicture}(-1,.5)(3,4)
\psline(1,0)(1,2)
\psline(1,2)(0,3)
\psline(1,2)(2,3)
\end{pspicture}}

\catcode`\@=11

\thinlines
\newskip\Einheit \Einheit=0.6cm
\newcount\xcoord \newcount\ycoord
\newdimen\xdim \newdimen\ydim \newdimen\PfadD@cke \newdimen\Pfadd@cke
\def\PfadDicke#1{\PfadD@cke#1 \divide\PfadD@cke by2 \Pfadd@cke\PfadD@cke \multiply\PfadD@cke by2}
\long\def\LOOP#1\REPEAT{\def\BODY{#1}\ITERATE}
\def\ITERATE{\BODY \let\next\ITERATE \else\let\next\relax\fi \next}
\let\REPEAT=\fi
\def\Punkt{\hbox{\raise-2pt\hbox to0pt{\hss\scriptsize$\bullet$\hss}}}
\def\DuennPunkt(#1,#2){\unskip
  \raise#2 \Einheit\hbox to0pt{\hskip#1 \Einheit
          \raise-2.5pt\hbox to0pt{\hss\normalsize$\bullet$\hss}\hss}}
\def\NormalPunkt(#1,#2){\unskip
  \raise#2 \Einheit\hbox to0pt{\hskip#1 \Einheit
          \raise-3pt\hbox to0pt{\hss\large$\bullet$\hss}\hss}}
\def\DickPunkt(#1,#2){\unskip
  \raise#2 \Einheit\hbox to0pt{\hskip#1 \Einheit
          \raise-4pt\hbox to0pt{\hss\Large$\bullet$\hss}\hss}}
\def\Kreis(#1,#2){\unskip
  \raise#2 \Einheit\hbox to0pt{\hskip#1 \Einheit
          \raise-4pt\hbox to0pt{\hss\Large$\circ$\hss}\hss}}
\def\Diagonale(#1,#2)#3{\unskip\leavevmode
  \xcoord#1\relax \ycoord#2\relax
      \raise\ycoord \Einheit\hbox to0pt{\hskip\xcoord \Einheit
         \unitlength\Einheit
         \line(1,1){#3}\hss}}
\def\AntiDiagonale(#1,#2)#3{\unskip\leavevmode
  \xcoord#1\relax \ycoord#2\relax \advance\xcoord by -0.05\relax
      \raise\ycoord \Einheit\hbox to0pt{\hskip\xcoord \Einheit
         \unitlength\Einheit
         \line(1,-1){#3}\hss}}
\def\Pfad(#1,#2),#3\endPfad{\unskip\leavevmode
  \xcoord#1 \ycoord#2 \thicklines\ZeichnePfad#3\endPfad\thinlines}
\def\ZeichnePfad#1{\ifx#1\endPfad\let\next\relax
  \else\let\next\ZeichnePfad
    \ifnum#1=1
      \raise\ycoord \Einheit\hbox to0pt{\hskip\xcoord \Einheit
         \vrule height\Pfadd@cke width1 \Einheit depth\Pfadd@cke\hss}%
      \advance\xcoord by 1
    \else\ifnum#1=2
      \raise\ycoord \Einheit\hbox to0pt{\hskip\xcoord \Einheit
        \hbox{\hskip-\PfadD@cke\vrule height1 \Einheit width\PfadD@cke depth0pt}\hss}%
      \advance\ycoord by 1
    \else\ifnum#1=3
      \raise\ycoord \Einheit\hbox to0pt{\hskip\xcoord \Einheit
         \unitlength\Einheit
         \line(1,1){1}\hss}
      \advance\xcoord by 1
      \advance\ycoord by 1
    \else\ifnum#1=4
      \raise\ycoord \Einheit\hbox to0pt{\hskip\xcoord \Einheit
         \unitlength\Einheit
         \line(1,-1){1}\hss}
      \advance\xcoord by 1
      \advance\ycoord by -1
    \else\ifnum#1=5
      \advance\xcoord by -1
      \raise\ycoord \Einheit\hbox to0pt{\hskip\xcoord \Einheit
         \vrule height\Pfadd@cke width1 \Einheit depth\Pfadd@cke\hss}%
    \else\ifnum#1=6
      \advance\ycoord by -1
      \raise\ycoord \Einheit\hbox to0pt{\hskip\xcoord \Einheit
        \hbox{\hskip-\PfadD@cke\vrule height1 \Einheit width\PfadD@cke depth0pt}\hss}%
    \else\ifnum#1=7
      \advance\xcoord by -1
      \advance\ycoord by -1
      \raise\ycoord \Einheit\hbox to0pt{\hskip\xcoord \Einheit
         \unitlength\Einheit
         \line(1,1){1}\hss}
    \else\ifnum#1=8
      \advance\xcoord by -1
      \advance\ycoord by +1
      \raise\ycoord \Einheit\hbox to0pt{\hskip\xcoord \Einheit
         \unitlength\Einheit
         \line(1,-1){1}\hss}
    \fi\fi\fi\fi
    \fi\fi\fi\fi
  \fi\next}
\def\hSSchritt{\leavevmode\raise-.4pt\hbox to0pt{\hss.\hss}\hskip.2\Einheit
  \raise-.4pt\hbox to0pt{\hss.\hss}\hskip.2\Einheit
  \raise-.4pt\hbox to0pt{\hss.\hss}\hskip.2\Einheit
  \raise-.4pt\hbox to0pt{\hss.\hss}\hskip.2\Einheit
  \raise-.4pt\hbox to0pt{\hss.\hss}\hskip.2\Einheit}
\def\vSSchritt{\vbox{\baselineskip.2\Einheit\lineskiplimit0pt
\hbox{.}\hbox{.}\hbox{.}\hbox{.}\hbox{.}}}
\def\DSSchritt{\leavevmode\raise-.4pt\hbox to0pt{%
  \hbox to0pt{\hss.\hss}\hskip.2\Einheit
  \raise.2\Einheit\hbox to0pt{\hss.\hss}\hskip.2\Einheit
  \raise.4\Einheit\hbox to0pt{\hss.\hss}\hskip.2\Einheit
  \raise.6\Einheit\hbox to0pt{\hss.\hss}\hskip.2\Einheit
  \raise.8\Einheit\hbox to0pt{\hss.\hss}\hss}}
\def\dSSchritt{\leavevmode\raise-.4pt\hbox to0pt{%
  \hbox to0pt{\hss.\hss}\hskip.2\Einheit
  \raise-.2\Einheit\hbox to0pt{\hss.\hss}\hskip.2\Einheit
  \raise-.4\Einheit\hbox to0pt{\hss.\hss}\hskip.2\Einheit
  \raise-.6\Einheit\hbox to0pt{\hss.\hss}\hskip.2\Einheit
  \raise-.8\Einheit\hbox to0pt{\hss.\hss}\hss}}
\def\SPfad(#1,#2),#3\endSPfad{\unskip\leavevmode
  \xcoord#1 \ycoord#2 \ZeichneSPfad#3\endSPfad}
\def\ZeichneSPfad#1{\ifx#1\endSPfad\let\next\relax
  \else\let\next\ZeichneSPfad
    \ifnum#1=1
      \raise\ycoord \Einheit\hbox to0pt{\hskip\xcoord \Einheit
         \hSSchritt\hss}%
      \advance\xcoord by 1
    \else\ifnum#1=2
      \raise\ycoord \Einheit\hbox to0pt{\hskip\xcoord \Einheit
        \hbox{\hskip-2pt \vSSchritt}\hss}%
      \advance\ycoord by 1
    \else\ifnum#1=3
      \raise\ycoord \Einheit\hbox to0pt{\hskip\xcoord \Einheit
         \DSSchritt\hss}
      \advance\xcoord by 1
      \advance\ycoord by 1
    \else\ifnum#1=4
      \raise\ycoord \Einheit\hbox to0pt{\hskip\xcoord \Einheit
         \dSSchritt\hss}
      \advance\xcoord by 1
      \advance\ycoord by -1
    \else\ifnum#1=5
      \advance\xcoord by -1
      \raise\ycoord \Einheit\hbox to0pt{\hskip\xcoord \Einheit
         \hSSchritt\hss}%
    \else\ifnum#1=6
      \advance\ycoord by -1
      \raise\ycoord \Einheit\hbox to0pt{\hskip\xcoord \Einheit
        \hbox{\hskip-2pt \vSSchritt}\hss}%
    \else\ifnum#1=7
      \advance\xcoord by -1
      \advance\ycoord by -1
      \raise\ycoord \Einheit\hbox to0pt{\hskip\xcoord \Einheit
         \DSSchritt\hss}
    \else\ifnum#1=8
      \advance\xcoord by -1
      \advance\ycoord by 1
      \raise\ycoord \Einheit\hbox to0pt{\hskip\xcoord \Einheit
         \dSSchritt\hss}
    \fi\fi\fi\fi
    \fi\fi\fi\fi
  \fi\next}
\def\Koordinatenachsen(#1,#2){\unskip
 \hbox to0pt{\hskip-.5pt\vrule height#2 \Einheit width.5pt depth1 \Einheit}%
 \hbox to0pt{\hskip-1 \Einheit \xcoord#1 \advance\xcoord by1
    \vrule height0.25pt width\xcoord \Einheit depth0.25pt\hss}}
\def\Koordinatenachsen(#1,#2)(#3,#4){\unskip
 \hbox to0pt{\hskip-.5pt \ycoord-#4 \advance\ycoord by1
    \vrule height#2 \Einheit width.5pt depth\ycoord \Einheit}%
 \hbox to0pt{\hskip-1 \Einheit \hskip#3\Einheit 
    \xcoord#1 \advance\xcoord by1 \advance\xcoord by-#3 
    \vrule height0.25pt width\xcoord \Einheit depth0.25pt\hss}}
\def\Gitter(#1,#2){\unskip \xcoord0 \ycoord0 \leavevmode
  \LOOP\ifnum\ycoord<#2
    \loop\ifnum\xcoord<#1
      \raise\ycoord \Einheit\hbox to0pt{\hskip\xcoord \Einheit\Punkt\hss}%
      \advance\xcoord by1
    \repeat
    \xcoord0
    \advance\ycoord by1
  \REPEAT}
\def\Gitter(#1,#2)(#3,#4){\unskip \xcoord#3 \ycoord#4 \leavevmode
  \LOOP\ifnum\ycoord<#2
    \loop\ifnum\xcoord<#1
      \raise\ycoord \Einheit\hbox to0pt{\hskip\xcoord \Einheit\Punkt\hss}%
      \advance\xcoord by1
    \repeat
    \xcoord#3
    \advance\ycoord by1
  \REPEAT}
\def\Label#1#2(#3,#4){\unskip \xdim#3 \Einheit \ydim#4 \Einheit
  \def\lo{\advance\xdim by-.5 \Einheit \advance\ydim by.5 \Einheit}%
  \def\llo{\advance\xdim by-.25cm \advance\ydim by.5 \Einheit}%
  \def\loo{\advance\xdim by-.5 \Einheit \advance\ydim by.25cm}%
  \def\o{\advance\ydim by.25cm}%
  \def\ro{\advance\xdim by.5 \Einheit \advance\ydim by.5 \Einheit}%
  \def\rro{\advance\xdim by.25cm \advance\ydim by.5 \Einheit}%
  \def\roo{\advance\xdim by.5 \Einheit \advance\ydim by.25cm}%
  \def\l{\advance\xdim by-.30cm}%
  \def\r{\advance\xdim by.30cm}%
  \def\lu{\advance\xdim by-.5 \Einheit \advance\ydim by-.6 \Einheit}%
  \def\llu{\advance\xdim by-.25cm \advance\ydim by-.6 \Einheit}%
  \def\luu{\advance\xdim by-.5 \Einheit \advance\ydim by-.30cm}%
  \def\u{\advance\ydim by-.30cm}%
  \def\ru{\advance\xdim by.5 \Einheit \advance\ydim by-.6 \Einheit}%
  \def\rru{\advance\xdim by.25cm \advance\ydim by-.6 \Einheit}%
  \def\ruu{\advance\xdim by.5 \Einheit \advance\ydim by-.30cm}%
  #1\raise\ydim\hbox to0pt{\hskip\xdim
     \vbox to0pt{\vss\hbox to0pt{\hss$#2$\hss}\vss}\hss}%
}
\catcode`\@=12

\begin{document}
 
\title{Non-commutative Hopf algebra of formal diffeomorphisms}
 
\author{Christian Brouder\footnote{\small \tt brouder@lmcp.jussieu.fr} \\
Laboratoire de Min\'eralogie-Cristallographie, \\
CNRS UMR7590, Universit\'es Paris 6 et 7, \\ 
IPGP, 4 place Jussieu, 75252 Paris Cedex 05, France; \\ 
{\em and\ } BESSY GmbH, Albert-Einstein-Str. 15,  
12489 Berlin, Germany. 
\and
Alessandra Frabetti\footnote{\small \tt frabetti@igd.univ-lyon1.fr} 
\ and
Christian Krattenthaler\footnote{\small \tt kratt@euler.univ-lyon1.fr}\ %
\footnote{Research partially supported by the
EC's IHRP Programme, grant HPRN-CT-2001-00272, 
``Algebraic Combinatorics in Europe."} \\
Institut Girard Desargues, CNRS UMR 5028, Universit\'e de Lyon 1,\\
B\^at. Braconnier, 26 av. Claude Bernard, 69622 Villeurbanne Cedex, 
France.}

\date{\today}
 
\maketitle
 
 
\begin{abstract}
The subject of this paper are two Hopf algebras which are the non-commutative 
analogues of two different groups of formal power series. The first group 
is the set of invertible series with the group law being multiplication 
of series, while the second group is the set of formal diffeomorphisms 
with the group law being composition of series. 
The motivation to introduce these Hopf algebras comes from the study of 
formal series with non-commutative coefficients. Invertible series with 
non-commutative coefficients still form a group, and we interpret the 
corresponding new non-commutative Hopf algebra as an alternative to the 
natural Hopf algebra given by the co-ordinate ring of the group, which has 
the advantage of being functorial in the algebra of coefficients. 
For the formal diffeomorphisms with non-commutative coefficients, this 
interpretation fails, because in this case the composition is not 
associative anymore. However, we show that for the dual non-commutative 
algebra there exists a natural co-associative co-product defining a 
non-commutative Hopf algebra. Moreover, we give an explicit formula 
for the antipode, which represents a non-commutative version of the 
Lagrange inversion formula, and we show that its coefficients are related 
to planar binary trees. 
Then we extend these results to the semi-direct co-product of the 
previous Hopf algebras, and to series in several variables. 
Finally, we show how the non-commutative Hopf algebras of formal series  
are related to some renormalization Hopf algebras, 
which are combinatorial Hopf algebras motivated by the renormalization 
procedure in quantum field theory, and to the renormalization functor 
given by the double tensor algebra on a bi-algebra.  
\end{abstract}
 
\newpage
\tableofcontents
\bigskip\bigskip

\section*{Introduction}
\addcontentsline{toc}{section}{\bf Introduction}

In the well-known paper \cite{ConnesKreimer}, A.~Connes and D.~Kreimer 
introduced a Hopf algebra structure on the set of Feynman graphs which 
allows one to describe the 
combinatorial part of the renormalization of quantum fields
in a very elegant and simple way. 
The renormalization procedure affects the coefficients of the perturbative 
series which describe the propagators and the coupling constants in 
quantum field theory, transforming them from infinite to well-defined finite 
quantities. 
The perturbative series involved are traditionally expanded over the 
set of Feynman graphs, but recent works by two of the authors 
\cite{Brouder,BFqedren}, showed that the amplitudes 
of Feynman graphs can be regrouped to form amplitudes associated to other 
combinatorial objects, such as rooted planar binary trees, 
or, at the coarsest level, 
positive integers. These new amplitudes correspond to 
new expansions of the perturbative series, and they turned out to be 
compatible with the renormalization. The renormalization is then encoded 
in the co-product of some Hopf algebras constructed on the set 
of rooted planar binary trees~\cite{BFqedtree} or on 
the set of positive integers~\cite{BF}.
The refinement of precision in the computation of the coefficients 
of the perturbative series, which is typical in quantum field theory, 
corresponds to a sequence of inclusions of Hopf algebras, the smallest one 
having generators 
labelled by the integers, the intermediate one 
with generators labelled by the trees, 
and the largest with generators labelled by the Feynman graphs. 

In this context, the use of Hopf algebras can be explained as a 
``local co-ordinate'' approach to the study of the renormalization groups, 
which are given on the sets of perturbative series relevant to each 
specific field theory. 
Accordingly, the Hopf algebras are never co-commutative, and they happen to 
be commutative if the amplitudes of the Feynman graphs are complex numbers, 
that is, if the quantum field considered is scalar. 
In this case, the renormalization Hopf algebras are exactly the co-ordinate 
rings of the renormalization groups. If the perturbative series are expanded 
over the integers, the co-products are constructed from the  
operations which are exactly the duals of the multiplication and 
the composition of usual formal series. 
On the contrary, when the series are expanded over trees, or over Feynman 
graphs, the dual operations of multiplication and composition 
of series have to be defined ``ad hoc'' in a way which generalises 
the usual operations. These new groups of series 
will be described in the upcoming paper \cite{Frabetti}
by one of the authors. 

However, there are cases in which it might be useful to consider series 
with non-scalar coefficients. This is the case, for instance, in 
quantum electro-dynamics \cite{BFqedren}, where the quantum fields are 
4-vectors or spinors, and the propagators are $4 \times 4$ matrices.
This is also the case for certain types of infrared renormalization 
\cite{BaganI,BaganII}, for the mass renormalization of a fermion family 
\cite{Bordes}, or for quantum field theory over noncommutative geometries 
\cite{PaschkeVerch}. 
For these latter cases, we are led to study the multiplication and 
the composition of formal series with non-scalar coefficients. 

In this paper, we consider two sets of formal power series with 
non-commutative coefficients. The first one is the set of invertible series 
with the multiplication law. Even with non-commutative coefficients, these 
series still form a group, and we show that its usual commutative 
co-ordinate ring can be replaced by a non-commutative Hopf algebra which 
is functorial in the algebra of coefficents. In fact, what we obtain is 
an example of a co-group element in the category of associative algebras, 
studied by B.~Fresse in~\cite{Fresse}, if we read it in the appropriate 
way, i.e., if we replace tensor products by free products among the algebras 
in the image of the co-product. 

The second set of formal series is that of formal diffeomorphisms on a line, 
with group law given by the composition of series. While this set 
forms a group if the coefficients of the series are scalar numbers, 
the composition fails to be associative when the coefficients are taken 
in an arbitrary non-commutative algebra. However, we show that on the dual 
algebra of local co-ordinates there is a natural co-product which is 
co-associative, and gives rise to a Hopf algebra which is neither commutative 
nor co-commutative. This Hopf algebra is related to the renormalization 
of quantum electrodynamics. More precisely, in Section~\ref{QEDHopf} we show 
that the non-commutative Hopf algebra of formal diffeomorphisms is a 
Hopf sub-algebra of the non-commutative Hopf algebra on planar binary trees,  
introduced in~\cite{BFqedtree}, which represents at the same time 
the charge and the photon renormalization Hopf algebras. 

Unlike the non-commutative Hopf algebra of invertible series, the 
non-commutative Hopf algebra of formal diffeomorphisms is not an example 
of a co-group element in associative algebras, because the co-product 
with image in the free product of algebras is not co-associative. 
However it has some remarkable properties, among which self-duality, 
which make it an oustanding example for many theories developed recently. 

For instance, in \cite{GavariniHdif}, F.~Gavarini applies the Quantum 
Duality Principle developed in \cite{GavariniQDP} to the non-commutative 
Hopf algebra of formal diffeomorphisms, via four one-parameter deformations, 
to get four quantum groups with semiclassical limits given by some Poisson 
geometrical symmetries. 

Similarly, P.~van der Laan describes in \cite{VanDerLaan} a general 
procedure to obtain canonically a Hopf algebra from an operad, and he shows 
that one obtains the non-commutative Hopf algebra of formal diffeomorphisms 
in the case of the operad of associative algebras. 

Finally, the non-commutative Hopf algebra of formal diffeomorphisms is 
the simplest example of a Hopf algebra obtained via the renormalization 
functor constructed on the double tensor algebra of a bi-algebra by 
W.~Schmitt and one of the authors in \cite{BrouderSchmitt}. This example 
will be treated in detail in Section~\ref{Rfunctor}. 
Because of this, the dual Hopf algebra of the Hopf algebra 
of formal diffeomorphisms 
is also the simplest non-trivial example of a Hopf algebra coming from a 
dendriform algebra, as introduced by J.--L.~Loday in \cite{LodayDialgebras}, 
and further studied with M.~Ronco in 
\cite{LodayRonco1,LodayRonco2,Ronco1,Ronco2}. 
In particular, the primitive elements of the dual Hopf algebra of 
formal diffeomorphisms are endowed with a very simple structure of a brace 
algebra (cf.~\cite{Ronco1,Ronco2}), which is under investigation by some 
of the authors. 

Our paper is organized as follows. 
In Section~\ref{series}, we recall the definition of the group of 
invertible series with the multiplication law, and of its co-ordinate ring. 
We first treat the case of series with scalar coefficients, and then 
consider series with non-commutative coefficients. We show that the 
classical co-ordinate ring can be replaced by a non-commutative Hopf algebra, 
and that the group can still be reconstructed as a group of characters 
with non-commutative values. 

In Section~\ref{diffeo}, we recall the definition of the group of 
formal diffeomorphisms with the composition law, and of its co-ordinate ring. 
Next, we consider series with non-commutative coefficients, and we show that, 
even if they do not form a group, on the algebra of local co-ordinates 
there is a natural structure of a Hopf algebra, which is neither 
commutative nor co-commutative. For this Hopf algebra, we present an explicit 
non-recursive formula for the antipode, which generalises the Lagrange 
Inversion Formula of formal series to the non-commutative context. 
We remark, that there appears already 
a non-commutative version of the Lagrange Inversion Formula in the
literature, which is due to I.~M.~Gessel \cite{GessAB}. However, the inversion
problem which is solved in \cite{GessAB} is inequivalent to ours.
Finally, we explain in detail the 
labelling of the antipode coefficients 
in terms of planar binary trees. 

In Section~\ref{semi-direct}, we show that the classical action of the 
group of formal diffeomorphisms on the group of invertible series 
can be generalized to the dual non-commutative context by means of a 
suitable co-action among Hopf algebras. This co-action allows one to construct 
the semi-direct (or smash) Hopf algebra of the previous Hopf algebras, 
studied by R.~K.~Molnar in~\cite{Molnar} and S.~Majid in~\cite{Majid}. 

The operations introduced in Section~\ref{semi-direct} are the ingredients 
which we need for explaining
the relationship between the non-commutative Hopf algebras 
of formal series and the renormalization Hopf algebras on planar binary trees 
used in~\cite{BFqedren,BFqedtree}. The latter Hopf algebras are
related to the renormalization 
of quantum electrodynamics. In Section~\ref{QEDHopf}, we prove that the 
algebras of series are Hopf sub-algebras of the corresponding algebras 
on trees. 

In Section~\ref{Rfunctor}, we show that the non-commutative Hopf algebra 
of formal diffeomorphisms can be obtained also via the renormalization 
functor described in~\cite{BrouderSchmitt}, applied to the simplest 
possible bi-algebra, the trivial one. On the one hand, this result places 
the non-commutative Hopf algebra of formal diffeomorphisms in the
context of a different 
approach of the renormalization procedure of quantum fields, namely the 
Epstein--Glaser renormalization on the configuration space, 
cf.~\cite{EpsteinGlaser}, and to its interesting further developement 
by G.~Pinter, cf.~\cite{Pinter1,Pinter2}. On the other hand, it relates the 
non-commutative Hopf algebra of formal diffeomorphisms to a large class 
of special algebras, such as dipterous, dendriform, brace and $B_\infty$ 
algebras, which were recently discovered by J.--L.~Loday and collaborators, 
cf.~\cite{LodayDialgebras,LodayRonco1,LodayRonco2,Ronco1,Ronco2}. 

Finally, in Section~\ref{variables}, we briefly sketch how to generalise 
the non-commutative Hopf algebra of formal diffeomorphisms to series 
with several variables. The practical applications of such formulae 
can be found in the renormalization of massive quantum field theory, 
cf.~\cite{BFmassren}.
 

\subsection*{Acknowledgments} 

The first two authors are very grateful to William Schmitt for several 
discussions on Hopf algebras and their antipodes. 
A.~Frabetti warmly thanks the Swiss National Foundation for Scientific 
Research for the support of her visit to the
Mathematics Departement of Lausanne 
University, where this work originated. She thanks as well 
the members of the Mathematics Departement of Lausanne 
University for their hospitality. 
Ch.~Brouder warmly thanks Alex Erko and the BESSY staff for their
hospitality.


\section{Non-commutative Hopf algebra of invertible series} 
\label{series}

In this section we introduce the Abelian group of invertible formal power 
series with scalar coefficients and its co-ordinate ring, which carries 
the structure of a commutative and co-commutative Hopf algebra. 
We use this simple example to describe the duality 
between a group and its co-ordinate ring. 

Then we introduce the group of invertible formal power series with 
coefficients from an arbitrary associative and unital algebra $\A$. 
Aside from the usual commutative co-ordinate ring, which depends on the 
chosen algebra $\A$, we present another Hopf algebra related to this 
group, which is no longer commutative and no longer dual to the group, 
but turns out to be functorial in $\A$.

\subsection{Group of invertible series}
\label{invertible}

Consider the set 
\begin{equation}
\label{Gp}
\Gi =\left\{ f(x) = 1 + \sum_{n=1}^\infty f_n x^n, \ f_n \in \C \right\}
\end{equation}
of invertible formal power 
series in a variable $x$ with complex coefficients, 
where, for simplicity, we fix the invertible constant term $f_0$ to be $1$. 
This set forms an Abelian group, with the multiplication 
\begin{equation}
\label{product}
(fg)(x) := f(x) g(x) 
= 1 + \sum_{n=1}^\infty \sum_{k=0}^n f_k g_{n-k}\ x^n, 
\end{equation}
the unit being given by the constant series $1(x)=1$, and where 
the inverse $f^{-1}$ of a series $f$ can be found recursively, 
for instance by using the Wronski formula (cf.\ \cite[ p.~17]{Henrici}). 
The first few coefficients of $f^{-1}$ are
$(f^{-1})_0 =1$, $(f^{-1})_1=-f_1$, $(f^{-1})_2=-f_2+f_1^2$. 

The group $\Gi$ is a projective limit of affine groups, so we can 
consider its co-ordinate ring $\C(\Gi)$ (defined below),
a commutative algebra 
related to the scalar functions on $\Gi$. 
We remark that, since $\Gi$ is not compact, 
and also not locally compact, its co-ordinate 
ring cannot be defined in the classical way as the algebra of 
representative functions on the group (i.e., the polynomials
in the matrix elements of the finite dimensional 
representations of the group, cf.\ \cite[Sec.~2.2]{Abe}). 
However, we can define $\C(\Gi)$ to be the set of functions $\Gi\to\C$
which are polynomial with respect to an appropriate basis.
As basis, we choose the functions $b_n$, $n=1,2,\dots$, where 
$b_n$ associates to each element of $\Gi$ its $n$-th coefficient. That
is, we may interpret $b_n$ as
the normalized $n^{th}$-derivative evaluated at $x=0$, 
\begin{equation*}
b_n (f) = \frac{1}{n!} \frac{\dd^n f(0)}{\dd x^n} = f_n. 
\end{equation*}
Thus, $\C(\Gi)$ is isomorphic to the polynomial ring
$\C[b_1,b_2,\dots]$.

The action of the functions on the elements of the group gives 
a duality pairing $\<b_n,f\>:= b_n(f)=f_n$ between $\C(\Gi)$ and $\Gi$. 
Through this pairing, the group structure on $\Gi$ induces the structure 
of a commutative Hopf algebra on $\C(\Gi)$, as it happens for affine or 
classical compact groups, cf.\ \cite{Abe,Hochschild,HewittRoss}. 
In particular, the group law on 
$\Gi$ induces a dual co-product $\Di$ on $\C(\Gi)$, that is, a map 
$\Di : \C(\Gi) \otimes \C(\Gi) \longrightarrow \C(\Gi)$ such that 
\begin{equation*}
\< \Di b_n , f \otimes g \> = \< b_n, fg \>,
\end{equation*}
where, as usual, $\<a\otimes b,f\otimes g\>=\<a,f\> \<b,g\>$.
In this case, it is easy to verify that the induced co-product on 
$\C(\Gi)$ has the form 
\begin{equation}
\label{Gp-coproduct}
\Di b_n = \sum_{k=0}^n b_k \otimes b_{n-k}, \quad (b_0:=1), 
\end{equation}
on the generators, and therefore it is co-commutative. 
Still by duality, the unit $1$ of the group induces a co-unit $\varepsilon$ 
on $\C(\Gi)$, that is, a map $\varepsilon : \C(\Gi) \longrightarrow \C$ 
such that 
\begin{equation*}
\varepsilon (b_n)= \< b_n,1 \>.  
\end{equation*}
Again, it is easy to verify that the co-unit has values 
$\varepsilon(1)=1$ and $\varepsilon(b_n)=0$ for $n\geq 1$. 
Finally, by duality, the operation of inversion in $\Gi$ gives rise to an 
antipode on $\C(\Gi)$, that is, a map $S: \C(\Gi) \longrightarrow \C(\Gi)$ 
such that 
\begin{equation*}
\< S (b_n),f \> = \< b_n , f^{-1} \>.  
\end{equation*}
In fact, the defining relation for the antipode yields a recursive
formula for the action of the antipode on the generators.
Thus, all these data together define the structure of a 
commutative and co-commutative graded connected Hopf algebra on $\C(\Gi)$. 

Finally, as is also the case for affine or classical Lie groups, cf.\ 
\cite{Krein,Tannaka} (see e.g.\ \cite[Ch.~VII, \S~30]{HewittRoss} or 
\cite[Theorem~3.5]{Hochschild}), the group $\Gi$ can be 
reconstructed completely from its co-ordinate ring $\C(\Gi)$, 
as the group of algebra homomorphisms (characters) 
$\Hom_{Alg}(\C(\Gi),\C)$ with the 
convolution product defined on the generators by 
\begin{equation}
\label{convolution}
        (\alpha \beta) (b_n) := m \circ(\alpha \otimes \beta) \circ\Di b_n,  
\end{equation}
for any algebra homomorphisms $\alpha,\beta$ 
on $\C(\Gi)$. Here, $m$ denotes the multiplication in $\C$. 
In other words, we have an isomorphism of groups 
$\Gi \cong \Hom_{Alg}(\C(\Gi),\C)$ which associates to a series 
$f \in \Gi$ the algebra homomorphism $\alpha_f$ on $\C(\Gi)$ given 
on the generators by $\alpha_f(b_n)= \<b_n,f\>$.


\subsection{Invertible series with non-commutative coefficients} 
\label{A-invertible}

Let $\A$ be an associative unital algebra, and consider the set 
\begin{equation}
\label{Gp(A)}
\Gi(\A) = \left\{ f(x) = 1 + \sum_{n=1}^\infty f_n x^n, 
\ f_n \in \A \right\}
\end{equation}
of invertible formal power series with coefficients in $\A$. 
The product \eqref{product} still makes $\Gi(\A)$ into a group, 
which is Abelian only if $\A$ is commutative.  

As before, $\Gi(\A)$ can be recovered from its co-ordinate ring 
$\C(\Gi(\A))$, which is still a polynomial ring  
in infinitely many variables, now depending on the chosen algebra $\A$. 
For instance, if $\A = M_2(\C)$ is the algebra of $2\times 2$ matrices 
with complex entries, then $f_n=\left(f_n^{ij}\right)_{i,j=1,2}$ 
with $f_n^{ij} \in \C$. Thus, we can choose the matrix elements as 
generators for the co-ordinate ring, and 
\begin{equation*}
\C(\Gi(M_2(\C))) \cong \C[b_1^{ij},b_2^{ij},...\ |\ i,j=1,2]
\end{equation*}
is a polynomial algebra on $4$ times infinitely many variables. 
The group law on $\Gi(\A)$ induces again a dual co-product on 
$\C(\Gi(\A))$. For instance, if $\A=M_2(\C)$ and 
$\C(\Gi(\A))=\C[b_n^{ij}|n\in\N,\ i,j=1,2]$, then 
\begin{equation*}
\D b_n^{ij}=\sum_{k=0}^n \sum_{l=1,2} b_k^{il} \otimes b_{n-k}^{lj}.
\end{equation*} 
As a result, $\C(\Gi(\A))$ is still a commutative Hopf algebra, 
which is co-commutative only if $\A$ is commutative. 

Alternatively, we can associate to $\Gi(\A)$ a non-commutative 
Hopf algebra $\Hi$, which has the advantage of being functorial 
in $\A$, that is, it does not depend on the chosen algebra $\A$. 
To do this, we consider the set $\Hi $ of 
$\A$-valued polynomial functions on $\Gi(\A)$. 
Then, as an algebra, $\Hi$ is isomorphic to the free unital associative 
algebra (tensor algebra) on infinitely many variables $b_n$, 
\begin{equation}
\label{Hp}
\Hi \cong \C\<b_1,b_2,\dots\>, 
\end{equation}
and the formula \eqref{Gp-coproduct} still defines a co-associative 
co-product which makes $\Hi$ into a non-commutative co-commutative 
Hopf algebra. Note that we can recover $\C(\Gi)$ from $\Hi$ by simple 
Abelianisation, and, thus, $\Hi$ can be considered as a non-commutative 
analogue of the group $\Gi$. 

Of course, in this case, the group $\Gi(\A)$ cannot anymore be 
reconstructed from $\Hi$, because the multiplication 
$m:\A \otimes \A \longrightarrow \A$ is not anymore an algebra homomorphism, 
and therefore the convolution product $\alpha \beta$ of two algebra 
homomorphisms $\alpha, \beta \in \Hom_{Alg}(\Hi,\A)$ defined by 
formula~\eqref{convolution} is not anymore an element of $\Hom_{Alg}(\Hi,\A)$. 
However, the group $\Gi(\A)$ can be easily reconstructed as follows. 

Let $\ast$ denote the free product of associative algebras 
(in the terminology of \cite{Loday:Kunneth} or \cite{Voiculescu}), 
which is the co-product or sum in the category of associative 
algebras (in the terminology of \cite[Section~I.7]{Lang}). 
We recall a few basic facts about the free product $\ast$. 

Given two unital associative algebras $\A$ and $\B$, the free product 
$\A\ast\B$ can be defined as the universal unital associative algebra 
which, for any given unital algebra $C$ and any algebra
homomorphisms $\alpha:\A\to\CC$ and $\beta:\B\to\CC$,
 makes the following diagram commutative:
$$
{\hbox{\setlength\unitlength{11pt}
\begin{picture}(13,6)
        \put(0,5){$\A$} 
        \put(1.2,5.2){\vector(1,0){3}}
        \put(1.2,4.6){\vector(1,-1){4}}
        \put(1.7,2){$\alpha$}  
        \put(4.9,5){$\A\ast\B$} 
        \put(6,4.5){\vector(0,-1){3}}
        \put(5,2.9){$\exists !$} 
        \put(5.8,0){$\CC$} 
        \put(11.5,5){$\B$} 
        \put(11,5.2){\vector(-1,0){3}}
        \put(11,4.6){\vector(-1,-1){4}}
        \put(9.9,2){$\beta$}
\end{picture}}} 
$$
As a vector space, $\A\ast\B$ is generated by the tensor products
in which elements of $\A$ and $\B$ alternate, that is
$$\A\ast\B = \bigg(\bigoplus_{k=0}^{\infty} 
(\A\otimes\B)^{\otimes k}\bigg)
\oplus\bigg(\B\otimes\bigoplus_{k=0}^{\infty} 
(\A\otimes\B)^{\otimes k}\bigg)
\oplus\bigg(\bigoplus_{k=1}^{\infty} 
(\B\otimes\A)^{\otimes k}\bigg)
\oplus\bigg(\bigoplus_{k=0}^{\infty} 
\A\otimes(\B\otimes\A)^{\otimes k}\bigg)
.$$ 
In particular, $\A\otimes\B$ is a sub-space of $\A\ast\B$. 
We endow $\A\ast\B$ with a product $\ast$, which is, essentially,
the concatenation product, except that any occurrence of 
$\cdots\otimes a\otimes a'\otimes\cdots$ is replaced by
$\cdots\otimes a a'\otimes\cdots$ for any $a,a'\in\A$, and
any occurrence of 
$\cdots\otimes b\otimes b'\otimes\cdots$ is replaced by
$\cdots\otimes b{} b'\otimes\cdots$ for any $b,b'\in\B$.
That is, for any $a,a'\in\A$ and $b,b'\in\B$, we have, for example,
$(a\otimes b) \ast (a'\otimes b')= a\otimes b\otimes a'\otimes b'$ and 
$(a\otimes b) \ast (b'\otimes a')= a\otimes (b b')\otimes a'$. 
Then, there is an obvious projection from $\A\ast\B$ to $\A\otimes\B$, 
which maps $a^1\otimes b^1\otimes a^2\otimes b^2\otimes \cdots
a^k\otimes b^k$ to $a^1a^2\cdots a^k\otimes b^1b^2\cdots b^k$, and
similarly for the other basis elements of $\A\ast \B$.
This map is a homomorphism of algebras. 

\begin{proposition}
\label{Hi*}
Let $\Dit:\Hi \longrightarrow \Hi \ast \Hi$ be the operator 
defined on the generators by formula~\eqref{Gp-coproduct}, and 
extended as an algebra homomorphism. Then $\Dit$ is co-associative. 

Moreover, if we denote by $\Hit$ the algebra $\Hi$ endowed with $\Dit$, 
then $\Hom_{Alg}(\Hit,\A)$ is a group with group law given by 
the convolution. In addition, the groups
$\Gi(\A)$ and $ \Hom_{Alg}(\Hit,\A)$ are isomorphic to each other. 
\end{proposition}

\begin{proof}
Since  $\Hi \otimes \Hi$ is a sub-space 
of $\Hi \ast \Hi$, formula~\eqref{Gp-coproduct} yields a well-defined 
operator on $\Hi$. 
The fact that $\Di$ is co-associative is easily checked, so it only remains 
to prove that $\Hom_{Alg}(\Hit,\A)$ is a group, and that it is
isomorphic to $\Gi(\A)$.  

The multiplication $m$ on $\A$, that is, the map $m:\A\otimes
\A\to\A$, can be extended to a map $m_\ast:\A\ast\A\to\A$. Unlike $m$,
the extension $m_\ast$ is an algebra homomorphism. Therefore, given
$\alpha,\beta\in \Hom_{Alg}(\Hit,\A)$, the convolution defined by
\begin{equation}
\label{*convolution}
\alpha \beta  :=  m_\ast \circ(\alpha \ast \beta) \circ\Di 
\end{equation}
is an element of $\Hom_{Alg}(\Hit,\A)$.
The rest of the proof that $\Hom_{Alg}(\Hit,\A)$ is a group, is
completely analogous to the proof in the commutative case, as, for
example, given in \cite{HewittRoss} or \cite{Hochschild}.

That $\Gi(\A)$ and $\Hom_{Alg}(\Hit,\A)$ are isomorphic to each other
is evident from the construction.
\end{proof}

Note that the new type of Hopf algebra $\Hit$ is an example 
of a co-group in the category of associative algebras, as considered 
by B.~Fresse \cite{Fresse} and by G.~M.~Bergman
and A.~O.~Hausknecht \cite[Sec.~60--62]{berghaus}. In particular, there
the reader may find more details on the group structure of
$\Hom_{Alg}(\Hit,\A)$ and on the generalization of this construction
to algebras over any operad.

Finally, note also that the co-product $\Di$ is just the composition
of $\Dit$ by the natural projection $\Hi\ast\Hi\to
\Hi\otimes\Hi$. Therefore, the co-associativity of $\Di$ follows from
the co-associativity of $\Dit$, but the converse is not true.


\section{Non-commutative Hopf algebra of series with composition}
\label{diffeo}

In this section we introduce the group of formal power series with the 
composition law, which we call {\it formal diffeomorphisms}, and
its co-ordinate ring. 

Proceeding as in Section~\ref{series}, we consider subsequently series with 
non-scalar coefficients and show that, even if these series do not
anymore form a group, dually there exists a Hopf algebra which is 
neither commutative nor co-commutative, and which reproduces the 
co-ordinate ring of the group by Abelianisation. 

For this new non-commutative Hopf algebra, we give an explicit formula
for the co-product 
and an explicit non-recursive formula for the antipode. 

\subsection{Group of formal diffeomorphisms and the Fa\`a di Bruno bi-algebra}
\label{diffeomorphisms}

We consider now the set 
\begin{equation}
\label{Gc}
\Gd=\left\{\varphi(x)=x+\sum_{n=1}^\infty \varphi_n x^{n+1}, \ 
\varphi_n \in \C \right\} 
\end{equation}
of formal power series in a variable $x$ with complex coefficients, 
zero constant term, and invertible linear term $\varphi_0$, 
which we set equal to $1$ for simplicity. 
This set forms a (non-Abelian) group with composition law 
\begin{equation}
\label{composition}
(\varphi \circ \psi)(x) := \varphi\big(\psi(x)\big) 
= \psi(x) + \sum_{n=1}^\infty \varphi_n \psi(x)^{n+1}. 
\end{equation}
The unit is given by the series $\id(x) = x$, and the (compositional) 
inverse $\varphi^{[-1]}(x)$ of a series $\varphi(x)$ can be found
by use of the Lagrange inversion formula \cite{Lagrange} (see e.g.\ 
\cite{Henrici} or \cite[Theorem~5.4.2]{Stanley}). 
Such series are called {\em formal diffeomorphisms} (tangent 
to the identity). 

An explicit expression for the composition can be easily derived
directly from the definition \eqref{composition}. However, we shall
not need it here.
Instead, for later use, we propose an alternative expression of 
the composition of two series, in form of the
({\it formal}) residue (see \cite{Egorychev} for an exposition of
formal residue calculus, in the commutative setting, however).
Given a Laurent series $F(z)$ in $z$, we write $\<z^{-1}\>F(z)$ for
the formal residue of $F(z)$, that is, its coefficient of $z^{-1}$.
Using this notation, the composition of $\varphi$ and $\psi$ can be
written as
\begin{equation}
\label{integral}
(\varphi\circ \psi)(x) = \res \frac{\varphi(z)}{z-\psi(x)}.
\end{equation}
This is justified, if we interpret $(z-\psi(x))^{-1}$ as a formal
power series in $z^{-1}$, that is, using the expansion of the
geometric series,
\begin{equation*}
\frac{1}{z-\psi(x)} = \frac{1}{z} 
\sum_{n=0}^\infty \left( \frac{\psi(x)}{z} \right)^n 
= \sum_{n=0}^\infty \psi(x)^n z^{-n-1} ,
\end{equation*} 
because then
\begin{equation*}
\frac{\varphi(z)}{z-\psi(x)} 
=\sum_{m=0}^\infty \sum_{n=0}^\infty \varphi_m \psi(x)^n z^{m-n},  
\end{equation*} 
from which \eqref{integral} follows immediately.
In the sequel, we shall {\it always} adopt this convention.

As before, this group is the projective limit of affine groups, and it 
can be reconstructed from its co-ordinate ring $\C(\Gd)$.  
The latter can be defined as the polynomial ring $\C[a_1,a_2,\dots]$ 
in infinitely many variables $a_n$, with $n \in \N$, where 
$a_n$ is the function on $\Gd$ acting 
as the normalized $(n+1)^{st}$-derivative evaluated at $x=0$, that is 
\begin{equation*}
a_n (\varphi) = \frac{1}{(n+1)!} 
\frac{\dd^{n+1} \varphi(0)}{\dd x^{n+1}} = \varphi_n. 
\end{equation*}
The group structure of $\Gd$ induces a Hopf algebra
structure on $\C(\Gd)$. The co-product for the generators
of $\C(\Gd)$ can be extracted from the standard duality condition 
\begin{equation*}
\<\Dd a_n,\varphi\otimes \psi\> = a_n(\varphi\circ \psi), 
\end{equation*}
where $\< a_n,\varphi\> = a_n(\varphi)$ and 
$\<a_n \otimes a_m , \varphi \otimes \psi\> = a_n(\varphi) a_m(\psi)$.

\begin{remark}
A slight variation of this Hopf algebra is known as the {\em Fa\`a di Bruno
bi-algebra}, an algebra based on the original computations made by Fa\`a 
di Bruno in \cite{FaadiBruno} on the derivatives of the composition 
of two functions. The Fa\`a di Bruno bi-algebra is in fact the coordinate 
ring of the semigroup of formal series of the form 
$\varphi(x) =\sum_{n=1}^\infty \varphi_n \frac{x^n}{n!}$, 
with $\varphi_1$ not necessarily equal to $1$. 
Repeating the duality procedure described above, we can identify  
the Fa\`a di Bruno bi-algebra in its standard form (cf.\ \cite{JoniRota} 
or \cite[Section~5.1]{Majid}) with the graded polynomial ring 
$\BFdB = \C[u_1,u_2,\dots]$ in infinitely many variables, with the degree of
$u_n$ being defined by $n-1$. The co-product in $\BFdB$, 
dual to the composition, 
takes the form 
\begin{align}
\label{D(un)}
\Delta u_n &= \sum_{k=1}^n u_k \otimes 
\underset{\alpha_1+\alpha_2+\dots+\alpha_n=k}
{\sum_{\alpha_1+2\alpha_2+\dots+n\alpha_n=n}} 
\frac{n!}{\alpha_1!\,\alpha_2!\cdots\alpha_n!}
\frac{u_1^{\alpha_1} u_2^{\alpha_2}\cdots u_n^{\alpha_n}}
{1!^{\alpha_1}2!^{\alpha_2}\cdots n!^{\alpha_n}}
\end{align}
on the generators $u_n$,
and the co-unit is defined by $\varepsilon(u_n)=\delta_{n,0}$.
For instance, 
\begin{align*}
\Delta u_1 &= u_1 \otimes u_1, \\
\Delta u_2 &= u_1 \otimes u_2 + u_2 \otimes u_1^2, \\
\Delta u_3 &= u_1 \otimes u_3 + u_2 \otimes 3 u_1 u_2 
+u_3 \otimes u_1^3, \\
\Delta u_4 &= u_1 \otimes u_4 + u_2 \otimes 4 u_1 u_3 
+ u_2 \otimes 3 u_2^2 + u_3 \otimes 6 u_1^2 u_2 + u_4 \otimes u_1^4.
\end{align*}

An explicit expression of the co-product $\Dd$ (in $\Hd$) 
on the generators $a_n$ 
can be obtained by replacing $\Delta$ by $\Dd$ and
$u_n$ by $n!\,a_{n-1}$  in \eqref{D(un)},
and by setting $u_1=a_0=1$.
\end{remark}

\begin{remark}
The co-ordinate ring $\C(\Gd)$ with its induced Hopf algebra structure
appeared also as a particular example of
an incidence Hopf algebra in the article \cite[Ex.~14.2]{Schmitt}
by W.~R.~Schmitt.
\end{remark}

\begin{remark}
A convenient way to present the co-product $\Dd$ for all generators $a_n$ in
compact form is by means of the generating series 
\begin{equation}
\label{A(x)} 
A(x) = x + \sum_{n=1}^\infty a_n x^{n+1} 
=\sum_{n=0}^\infty a_n x^{n+1}, \quad (a_0 :=1),
\end{equation}
since it allows to reconstruct each series $\varphi \in \Gd$ by duality: 
\begin{equation} \label{calc1}
\< A(x) ,\varphi\>
= \sum_{n=0}^\infty  \< a_n,\varphi \> x^{n+1}
= \sum_{n=0}^\infty \varphi_n x^{n+1} = \varphi(x). 
\end{equation}
In fact, more generally, we have 
\begin{align} \notag
\< A(x)^m,\varphi \> &= 
\left \< \sum _{n_1,\dots,n_m\ge 0} ^{}a_{n_1}a_{n_2}\cdots
a_{n_m}x^{(n_1+1)+\dots+(n_m+1)} ,\varphi\right\> \\
\notag
&= \sum _{n_1,\dots,n_m\ge 0} ^{}  \< a_{n_1}a_{n_2}\cdots
a_{n_m}  ,\varphi\> x^{(n_1+1)+\dots+(n_m+1)}\\
&= \sum _{n_1,\dots,n_m\ge 0} ^{}  \varphi_{n_1}\varphi_{n_2}\cdots
\varphi_{n_m}  x^{(n_1+1)+\dots+(n_m+1)}
= \varphi(x)^m. \label{calc2}
\end{align}
If we now set $\Dd A(x) := \sum \Dd a_n x^n$, then we obtain 
\begin{align*}
\<\Dd A(x), \varphi\otimes \psi \> &=
\sum_{n=0}^\infty \< \Dd a_n,\varphi \otimes \psi \> x^{n+1}
= \sum_{n=0}^\infty a_n(\varphi\circ \psi) x^{n+1}
= (\varphi\circ \psi)(x) \\ 
&= \res \frac{\varphi(z)}{z-\psi(x)}\\
& = \res\left(\<A(z),\varphi\>\, \<\frac{1}{z-A(x)},\psi \>\right)
= \res \< A(z) \otimes \frac{1}{z-A(x)},\varphi\otimes\psi \> ,
\end{align*}
where we have used \eqref{calc1} and \eqref{calc2} to go from the second to
the third line.
Therefore, the co-product of the generating series $A(x)$ is given by 
\begin{equation}
\label{A(x)-coproduct}
\Dd A(x) = \res A(z) \otimes \frac{1}{z-A(x)} .
\end{equation}
To find $\Dd a_n$ for each $n \geq 1$, it suffices to evaluate 
this residue, where, again, the inverse of $z-A(x)$ has to be
interpreted as
\begin{equation*}
\frac{1}{z-A(x)} = \frac{1}{z} 
\sum_{n=0}^\infty \left( \frac{A(x)}{z} \right)^n 
= \sum_{n=0}^\infty A(x)^n z^{-n-1} . 
\end{equation*} 
\end{remark}


\subsection{Formal diffeomorphisms with non-commutative coefficients} 
\label{sec:formdiff}

Let $\A$ be an associative unital algebra, and consider the set 
\begin{equation*}
\label{Gc(A)}
\Gd(\A)=\left\{ \varphi(x)=\sum_{n=0}^\infty \varphi_n x^{n+1},
\varphi_n \in \A, \varphi_0=1 \right\} 
\end{equation*}
of formal power series with coefficients in $\A$.
Proceeding in analogy to Section~\ref{A-invertible}, where we adopted
formula \eqref{product} for invertible series for the case of not necessarily
commutative coefficients, it seems natural to adopt formula 
\eqref{composition} as the definition of the composition 
$\varphi \circ \psi$ for series $\varphi$ and $\psi$ in not necessarily 
commutative coefficients. However, such a composition is not associative 
unless $\A$ is commutative, since the associator 
\begin{equation*}
\big(\varphi\circ(\psi\circ\eta)\big)(x)
- \big((\varphi\circ\psi)\circ\eta)\big)(x)
= x^4 \big(\varphi_1\eta_1\psi_1-\varphi_1\psi_1\eta_1) + O(x^5)
\end{equation*}
is non-zero if the coefficients do not commute. 

Thus, the set of series with non-commutative coefficients, with
composition defined by \eqref{composition}, does not 
define a group. However, we shall see in this section that, dually, 
there exists a co-associative co-product which gives rise to a 
non-commutative Hopf algebra.
In other words, if we consider the free non-commutative algebra 
generated by the elementary functions $a_n(\varphi)=\varphi_n \in \A$, 
for $\varphi \in \Gd(\A)$, then the dual co-product of the composition 
\eqref{composition} is co-associative. 

\begin{definition} \label{Adef}
Let $\Hd = \C\<a_1,a_2,\dots\>$ denote the free associative algebra 
in infinitely many variables $a_n$, for $n \in \N$. 
As for the commutative case \eqref{diffeomorphisms}, we consider 
the generating series $A(x) = x + \sum_{n=1}^\infty a_n x^{n+1}$. 
Then we define a co-product on the generators of $\Hd$ by the global 
formula \eqref{A(x)-coproduct},
\begin{equation} \label{eq:DA}
\Dd A(x) = \res A(z) \otimes \frac{1}{z-A(x)} , 
\end{equation}
and we extend it multiplicatively to products of elements of $\Hd$. 
As a co-unit on $\Hd$, we take the standard graded co-unit 
$\varepsilon(1)=1$ and $\varepsilon(a_n)=0$, in other words, 
$\varepsilon(A(x))=x$. 
We show in Theorem~\ref{Hd-Hopf} that $\Hd$ is indeed a Hopf algebra. 
Of course, we can obtain $\C(\Gd)$ from $\Hd$ by taking the 
Abelianisation. 
\end{definition} 

For the structural analysis of $\Hd$ and its co-product $\Dd$, we
shall make frequent use of certain non-commutative polynomials, which
we define next.

\begin{definition} \label{defQ}
For $m,n\ge0$, we define the polynomials $Q_m^{(n)}(a)$ 
in $m$ variables $a_1,a_2,\dots,a_m$ by
\begin{equation*}
Q_m^{(n)}(a) = \underset{j_0,\dots,j_{n} \ge 0}{ \sum_{j_0+\cdots+j_{n}=m}}
a_{j_0} \cdots a_{j_n},
\end{equation*}
where, as before, $a_0$ is interpreted as $1$. For convenience, we set
$Q_m^{(-1)}(a)=0$ if $m>0$ and $Q_0^{(-1)}(a)=1$.
\end{definition} 

According to this definition, we have $Q^{(n)}_0(a)=1$ for all $n$,
$Q^{(0)}_m(a)=a_m$ for all $m$, and
\begin{align*}
Q^{(n)}_1(a) &= (n+1) a_1, \\ 
Q^{(n)}_2(a) &= (n+1) a_2 + \frac{n(n+1)}{2} a_1^2, \\ 
Q^{(n)}_3(a) &= (n+1) a_3 + \frac{n(n+1)}{2} (a_1a_2+a_2a_1) + 
\frac{(n+1)n(n-1)}{6} a_1^3, 
\end{align*}
for any $n \geq 0$. More generally, for $n\ge0$, 
we may write $Q_m^{(n)}(a)$ in the
form
\begin{equation}
\label{polP}
Q_m^{(n)}(a) =
\sum _{l=0} ^{\infty} \binom {n+1}{l}
\underset{h_1,\dots,h_{l} \ge 1}{ \sum_{h_1+\cdots+h_{l}=m}}
a_{h_1} \cdots a_{h_l}.
\end{equation}

It follows directly from the definition that
\begin{align} \label{eq:Qcoef}
A(x)^{n+1} = x^{n+1}\left( 1 +\sum_{p=1}^\infty a_p x^{p} \right)^{n+1} 
= \sum_{m=0}^\infty Q_m^{(n)}(a) \ x^{m+n+1} ,
\end{align}
in other words, $Q_m^{(n)}(a)$ is the coefficient of $x^{m+n+1}$ in
$A(x)^{n+1}$. From this generating function, we can easily derive two
equations satisfied by 
the $Q_m^{(n)}(a)$ which we shall frequently use later on.

\begin{lemma}
\label{recursiveQ}
For $n,m\ge0$, the polynomials $Q^{(n)}_{m}(a)$ satisfy the recurrence
\begin{equation*}
Q^{(n)}_{m} (a)=  \sum_{l=0}^{m} a_l Q^{(n-1)}_{m-l}(a) =
\sum_{l=0}^{m}  Q^{(n-1)}_{m-l}(a)a_l .
\end{equation*}
\end{lemma}
 
\begin{proof} 
The equation follows directly by comparing coefficients of $x^{n+m+1}$ in the
generating function identity $A(x)^{n+1}=A(x)A(x)^n=A(x)^nA(x)$, using
\eqref{eq:Qcoef}.
\end{proof}

\begin{lemma}
\label{quadraticQ}
For $l,m,n\geq 0$, the polynomials $Q^{(n)}_m(a)$ satisfy the quadratic 
relation  
\begin{align}
Q^{(l+n+1)}_m (a)&= \sum_{k=0}^{m} Q^{(l)}_k(a) Q^{(n)}_{m-k}(a).
 \label{lemma3} 
\end{align}
This relation holds as well if one of\/ $l$ and $n$ is equal to $-1$.
\end{lemma}
 
\begin{proof} 
The relation follows directly by comparing coefficients of
$x^{m+l+n+2}$ in the generating function equation
$A(x)^{l+n+2}=A(x)^{l+1}A(x)^{n+1}$, using \eqref{eq:Qcoef}.
\end{proof}
 
Another simple corollary of \eqref{eq:Qcoef} is an explicit expression
for the double generating function for the $Q^{(n)}_{m}(a)$.

\begin{corollary} 
\label{genQ(x,y)}
The generating function 
\begin{align}
Q(x,y) &= \sum_{m=0}^\infty\sum_{n=0}^\infty x^n y^m Q^{(n)}_m(a) 
\label{defG}
\end{align}
is the Green function\footnote{%
If $H$ is an operator on a Hilbert space, the {\em resolvent} or 
{\em Green function} of $H$ is the operator $R(z)=(H-z)^{-1}$, 
and it satisfies the {\em resolvent identity} $R(z)=R(z')+(z-z')R(z)R(z')$, 
cf.~\cite{Kato}.}  
at $x$ of the Hamiltonian $H(y)=\left(\frac {A(y)} {y}\right)^{-1}$, 
where, as before, $A(y)=\sum_{n=0}^\infty a_n y^{n+1}$.
That is, the generating function is 
\begin{align}
\label{GenGreen}
Q(x,y) &= \left(\left(\frac {A(y)} y\right)^{-1} - x\right)^{-1} 
= \left(1 - x \frac {A(y)} y\right)^{-1} \frac {A(y)} y . 
\end{align}
Moreover, it satisfies the resolvent equation
\begin{align}
\label{resolvent}
Q(x,y) &= Q(z,y) + (x-z) Q(x,y) Q(z,y). 
\end{align}
\end{corollary}
 
\begin{proof} 
For proving the first equation, we use \eqref{eq:Qcoef} to rewrite
$Q(x,y)$ in the form
\begin{align*}
Q(x,y) &=  \sum_{n=0}^\infty x^n \sum_{m=0}^\infty y^m 
 Q^{(n)}_{m}(a)=  \sum_{n=0}^\infty  x^n \left(
\frac {A(y)} {y}\right)^{n+1} .
\end{align*}
The sum over $n$ is a geometric series and can therefore be evaluated.
This yields \eqref{GenGreen}. Eq.~\eqref{resolvent} can now easily be
verified by substituting \eqref{GenGreen} in \eqref{resolvent}.
\end{proof}

We are now in the position to explicitly describe the action
of $\Dd$ on the generators $a_n$.

\begin{lemma} 
\label{Dd(a_n)}
On the generators $a_n$ of $\Hd$, the co-product is given by
\begin{align*}
\Dd a_n &=
\sum_{k=0}^{n} a_{k} \otimes Q^{(k)}_{n-k}(a).
\end{align*}
\end{lemma}

\begin{proof} 
According to Definition~\ref{Adef}, we find $\Dd A_n$ by extracting
the coefficient of $x^{n+1}$ from $\Dd A(x)$, as given by
\eqref{eq:DA}. Consequently, we expand the right-hand side of
\eqref{eq:DA}, and we obtain
\begin{align}
\Dd A(x) &= \res A(z) \otimes \frac{1}{z-A(x)} \notag \\ 
&= \sum_{n=0}^\infty \res z^{-n-1} A(z) \otimes A(x)^n \notag \\ 
&= \sum_{n=0}^\infty \sum_{k=0}^\infty 
a_k \otimes A(x)^n \res z^{k-n} \notag \\ 
&= \sum_{k=0}^\infty a_k \otimes A(x)^{k+1} \label{DdA(x)}.
\end{align}
Use of \eqref{eq:Qcoef} thus yields our claim.
\end{proof}

For instance, we have
\begin{align*}
\Dd a_1 &= a_1 \otimes 1 + 1 \otimes a_1, \\
\Dd a_2 &= a_2 \otimes 1 + 1 \otimes a_2 + 2 a_1 \otimes a_1 , \\
\Dd a_3 &= a_3 \otimes 1 + 1 \otimes a_3
+ 3 a_2 \otimes a_1 + 2 a_1 \otimes a_2 + a_1 \otimes a_1^2 .
\end{align*}
{}From the value of $\Dd a_3$, we see that the co-product
is not co-commutative.

There is, as well, an elegant description of the co-product on the
polynomials $Q^{(n)}_m(a)$, which we give in the corollary below. 
In particular, it shows that the sub-algebra of $\Hd$ generated by the
$Q^{(n)}_m(a)$'s is in fact a {\it Hopf} sub-algebra of $\Hd$.

\begin{corollary}
\label{DeltaQ}
For $m,n\ge0$, the co-product of $Q^{(n)}_m(a)$ is given by 
\begin{align} \label{DifQ}
\Dd Q^{(n)}_{m} (a)&= 
 \sum_{k=0}^{m} Q^{(n)}_{m-k}(a) \otimes Q^{(n+m-k)}_{k}(a).
\end{align}
\end{corollary}
 
\begin{proof}
Since, by \eqref{eq:Qcoef}, the generating function of the
$Q^{(n)}_m(a)$ is $A(x)^{n+1}$, we compute the image of a power of
$A(x)$ under $\Dd$. 
Using \eqref{eq:DA}, we have
\begin{align*}
\Dd \big(A(x)^2\big) &= \Dd A(x) \ \Dd A(x) 
= \sum_{m,n=0}^\infty a_m a_n \otimes A(x)^{m+n+2} \\ 
&= \res \sum_{m,n,k=0}^\infty 
a_m a_n z^{m+n+2}\otimes A(x)^k z^{-k-1} \\ 
&= \res A(z)^2 \otimes \frac{1}{z-A(x)}. 
\end{align*}
By induction on $n$, a similar reasoning for $n \geq 2$ shows that
\begin{equation} \label{DifA}
\Dd \big(A(x)^n\big) = 
\res A(z)^n \otimes \frac{1}{z-A(x)}. 
\end{equation}
Hence, comparison of coefficients of $x^{m+n+1}$ in \eqref{DifA}
(with $n$ replaced by $n+1$) yields
\eqref{DifQ}. 
\end{proof}

The results obtained so far, allow us not to prove that $\Hd$ is a
graded and connected Hopf algebra.

\begin{theorem} 
\label{Hd-Hopf}
The algebra $\Hd$ is a
graded and connected Hopf algebra, which is
 neither commutative nor co-commutative.
\end{theorem}

\begin{proof}
The algebra $\Hd$ becomes a graded algebra and a graded co-algebra 
by defining the degree of a
monomial 
$a_{j_1}a_{j_2}\cdots a_{j_m}$ to be $j_1+j_2+\dots +j_m$.
Moreover, $\Hd$ is connected, that is, the zero degree part consists
only of the scalars.
Since, for a graded
connected Hopf algebra, the antipode is given by 
the standard recursive formula 
\begin{equation} \label{recform}
\Sd a_n = -a_n - \sum_{p=1}^{n-1} a_p \Sd \big(Q^{(p)}_{n-p}(a)\big)  
= -a_n - \sum_{p=1}^{n-1} (\Sd a_p) Q^{(p)}_{n-p}(a) , 
\end{equation}
the only statement which needs to be proved is the co-associativity 
of the co-product. 

A recursive proof of the latter is given in \cite{BF}.
It involves a number of technical computations 
in order to deduce some recurrence relations 
for the polynomials $Q_m^{(l)}(a)$. 
Here we present an alternative simple proof which justifies 
the introduction of residues and generating series. 

If we consider a series $f(x)=\sum_{n=0}^\infty f_n x^n$ 
with complex coefficents $f_n\in\C$, then, because of \eqref{DifA}, 
the image of the the composition 
$(f \circ A)(x) =\sum_{n=0}^\infty f_n A(x)^n$ under $\Dd$ is given by
\begin{equation*}
\Dd (f\circ A)(x) = 
\res (f\circ A)(z) \otimes \frac{1}{z-A(x)}.
\end{equation*}
In particular, for $f(x)=1/(y-x)$ (regarded as a formal power series in $x$!)
we obtain 
\begin{equation}
\Dd \frac{1}{y-A(x)} =
\res \frac{1}{y-A(z)} \otimes \frac{1}{z-A(x)}. \label{Ddfrac}
\end{equation}

With the help of this result, it is now very easy to prove the 
co-associativity of the co-product: on the one hand, we have 
\begin{align*}
(\Dd\otimes 1)\circ \Dd A(x) &= 
\resz{z_1} \Dd A(z_1) \otimes \frac{1}{z_1-A(x)} \\ 
&= \resz{z_1} \resz{z_2} 
A(z_2) \otimes \frac{1}{z_2-A(z_1)} \otimes \frac{1}{z_1-A(x)}, 
\end{align*}
and, on the other hand, we have 
\begin{align*}
(1 \otimes \Dd)\circ\Dd A(x) &= 
\resz{z_2} A(z_2) \otimes \Dd\frac{1}{z_2-A(x)} \\ 
&= \resz{z_2} \resz{z_1} 
A(z_2) \otimes \frac{1}{z_2-A(z_1)} \otimes \frac{1}{z_1-A(x)},
\end{align*}
yielding the same expression on the right-hand side.
\end{proof}

\begin{remark}
The Abelianisation of $\Hd$ gives the co-ordinate ring $\C(\Gd)$. 
In fact, if we suppose that the variables $a_n$ commute, 
then, rewriting \eqref{polP}, the polynomial $Q_m^{(k)}(a)$ becomes 
\begin{equation*}
Q_m^{(k)}(a) = 
\sum _{l=0} ^{\infty} \binom {k+1}{l}
\underset
{p_1+p_2+\dots +p_m=l}{\sum_{p_1+2p_2+\dots +mp_m=m}}
\frac{l!}{p_1!\cdots p_m!} \ a_1^{p_1}\cdots a_m^{p_m},
\end{equation*}
where the sum runs over $p_1,\dots,p_m\geq 0$.\footnote{
It should be noted that these combinatorial factors are well known
from Planck's quantum theory of blackbody radiation, see for instance 
\cite{Planck}.} 
Consequently, the co-product $\Dd$, as given Lemma~\ref{Dd(a_n)}, 
becomes 
\begin{align}
\label{D(an)}
\Dd a_n &= 
\sum_{k=0}^{n} a_k \otimes \sum_{l=1}^{n-k} \frac{(k+1)!}{(k+1-l)!} 
\underset{p_1+p_2+\dots +p_{n-k}=l}{\sum_{p_1+\dots +(n-k)p_{n-k}=n-k}} 
\frac{1}{p_1!\cdots p_{n-k}!} \ a_1^{p_1}\cdots a_{n-k}^{p_{n-k}}
\end{align}
in the commutative case.
This agrees with the Fa\`a di Bruno co-product on the variables 
$a_n=\break u_{n+1}/(n+1)!$, with $a_0=u_1=1$. 
In fact, from \eqref{D(un)}, by setting $\beta_i=\alpha_{i+1}$, and then 
summing up over $l=\beta_0$, we obtain: 
\begin{align*}
\Delta a_n &= \sum_{k=0}^n a_k \otimes 
\underset{\beta_0+\beta_1+\dots+\beta_n=k+1}
{\sum_{\beta_0+2\beta_1+\dots+(n+1)\beta_n=n+1}} 
\frac{(k+1)!}{\beta_0!\beta_1!\cdots\beta_n!}
a_1^{\beta_1} a_2^{\beta_2}\cdots a_n^{\beta_n} \\ 
&= \sum_{k=0}^n a_k \otimes \sum_{l=0}^k \frac{(k+1)!}{l!} 
\underset{\beta_1+\dots+\beta_n=k+1-l}
{\sum_{2\beta_1+\dots+(n+1)\beta_n=n+1-l}} 
\frac{1}{\beta_1!\cdots\beta_n!}
a_1^{\beta_1} a_2^{\beta_2}\cdots a_n^{\beta_n} \\ 
&= \sum_{k=0}^n a_k \otimes \sum_{l=1}^{k+1} \frac{(k+1)!}{(k+1-l)!} 
\underset{\beta_1+\dots+\beta_n=l}
{\sum_{\beta_1+\dots+n\beta_n=n-k}} 
\frac{1}{\beta_1!\cdots\beta_n!}
a_1^{\beta_1} a_2^{\beta_2}\cdots a_n^{\beta_n}. 
\end{align*}
Since the $\beta_i$ are non-negative integers, the condition 
$\beta_1+\dots+n\beta_n=n-k$ implies that $\beta_{n-k+1}=\dots=\beta_n=0$. 
As a consequence, the condition $\beta_1+\dots+\beta_{n-k}=l$ implies that 
$k$ runs from $1$ to at most $n-k$. Therefore $\Delta a_n$ gives exactly 
\eqref{D(an)}.  
\end{remark}

\begin{remark}
In Proposition~\ref{Hi*}, we showed that the co-product $\Di$ on the 
Hopf algebra $\Hi$ can be lifted up to a new kind of co-product 
$\Dit$, with values in the free product $\Hi \ast \Hi$, which remains 
co-associative and gives $\Hit = (\Hi,\Dit)$ the structure of a 
co-group in associative algebras, cf.~\cite{Fresse}. 
By way of contrast, a lifting of the co-product $\Dd$ with values in 
$\Hd \ast\Hd$ is {\em not\/} co-associative, and the new Hopf algebra 
$\Hdt=(\Hd,\Ddt)$ is {\em not\/} a co-group in associative algebras. 
This reflects the fact that the set of formal diffeomorphisms with 
non-commutative coefficients fails to be a group because the composition 
fails to be associative. 
In fact, given a non-commutative algebra $\A$, the set $\Gd(\A)$ 
can still be reconstructed from $\Hdt$ as the set $\Hom_{Alg}(\Hdt,\A)$ 
of characters with non-commutative values, via the adapted convolution 
defined by formula~\eqref{*convolution}. In this duality, 
the non-co-associativity of $\Ddt$ corresponds exactly to the 
non-associativity of the composition of series. 
This makes it even more remarkable that the co-product $\Dd$, which 
is the composition of $\Ddt$ by the natural projection 
$\Hd\ast\Hd \to \Hd\otimes\Hd$, turns out to be co-associative. 

Finally note that the problem of whether the Hopf algebra on trees introduced 
by J.--L.~Loday and M.~Ronco in \cite{LodayRoncoHopfTree} is a co-group
was already investigated by R.~Holtkamp in \cite{Holtkamp}. 
This Hopf algebra turns out to be isomorphic 
to the linear dual of the Hopf algebra $\Hat$ that we introduce in 
Section~\ref{QEDHopf}, which we show being an extension of $\Hd$. 
In Sections~3.5 and 3.6 of \cite{Holtkamp}, Holtkamp shows that the 
Loday--Ronco Hopf algebra cannot be a co-group,
and therefore our result on $\Hd$ agrees with his. 
\end{remark}


\subsection{Explicit non-commutative Lagrange formula} 

In the commutative setting, the antipode of the co-ordinate ring 
$\C(\Gd)$ is of course the operation which is dual to the inversion 
$\phi \mapsto \phi^{[-1]}$ of a formal power series. 
The coefficients of the inverse series $\phi^{[-1]}$ are usually computed 
by using the Lagrange inversion formula.

As we outlined in Section~\ref{sec:formdiff},
in the non-commutative setting, the analogue for the co-ordinate $\C(\Gd)$
is the (non-commutative) Hopf algebra $\Hd$. Its antipode 
can be computed using the
recursive formula \eqref{recform}. 
In particular, in $\Hd$, the first few values of the antipode 
on the generators $a_i$ are
\begin{align*}
\Sd a_1 &= -a_1 , \\
\Sd a_2 &= -a_2 + 2 a_1^2 , \\
\Sd a_3 &= -a_3 + (2 a_1 a_2 + 3 a_2 a_1) - 5a_1^3, \\
\Sd a_4 &= -a_4 + (2 a_1 a_3 + 3 a_2^2 + 4 a_3 a_1) 
- (5 a_1^2 a_2 + 7 a_1 a_2 a_1 + 9 a_2 a_1^2) + 14 a_1^4. 
\end{align*}
It should be observed that the square of the antipode is equal to the identity 
only for $a_1$ and $a_2$. 

In contrast to the commutative setting, where there exist both the
recursive and the Lagrange inversion formula, for $\Hd$ there is no known 
analogue of the Lagrange inversion formula. In the theorem below, we
present an explicit non-recursive expression for the antipode for $\Hd$. 
Clearly, the Abelianisation of this expression gives an explicit formula
for the inversion of formal diffeomorphisms (which can be seen as an
alternative version of the Lagrange inversion formula). 

\begin{theorem}
The action of the antipode of the Hopf algebra $\Hd$  
on the generators has the following closed form: 
\begin{equation}
\label{antipode}
\Sd a_n = -a_n -\sum_{k=1}^{n-1} (-1)^k 
\underset {n_1,\dots,n_{k+1} >0}{  \sum_{n_1+\cdots+n_k+n_{k+1}=n}} 
\ll(n_1,\dots,n_k)\ a_{n_1} \cdots a_{n_k} a_{n_{k+1}} , 
\end{equation} 
with coefficients 
\begin{equation} 
\label{ll}
\ll(n_1,\dots,n_k) = 
\underset{h=1,\dots,k-1}{\underset{m_1+\dots+m_h\ge h}
{\sum_{m_1+\dots+m_k=k}}} 
\binom{n_1+1}{m_1} \cdots \binom{n_k+1}{m_k} .
\end{equation} 
\end{theorem}

It should be noted that the coefficient $\ll(n_1,\dots,n_k)$ of the monomial 
$a_{n_1} \cdots a_{n_k} a_{n_{k+1}}$ does not depend on the last factor 
$a_{n_{k+1}}$. 

\begin{proof}
We prove formula \eqref{antipode} by induction on $n$, 
starting from the recursive
definition \eqref{recform} of the antipode $\Sd$.
Formula \eqref{antipode} is obviously correct for $n=1$, and it is
also correct for
$n=2$, because in the latter case 
the sum over $k$ on the right-hand side gives 
only one term for $k=1$, namely 
\begin{equation*}
(-1) \underset{n_1,n_2 >0}{ \sum_{n_1+n_2=2}} \ll(n_1)\ a_{n_1} a_{n_2}
= - \ll(1)\ a_1^2 = - \binom{2}{1}\ a_1^2 = -2 a_1^2. 
\end{equation*}

For $n \geq 3$, suppose that formula \eqref{antipode} holds for 
$\Sd a_p$, $p=1,\dots,n-1$. We shall show that 
\begin{equation}
\label{halfantipode}
\sum_{p=1}^{n-1} (\Sd a_p) Q^{(p)}_{n-p} (a)= 
\sum_{p=1}^{n-1} (-1)^p 
\underset{n_1,\dots,n_{p+1} >0}{\sum_{n_1+\cdots+n_p+n_{p+1}=n}} 
\ll(n_1,\dots,n_p)\ a_{n_1} \cdots a_{n_p} a_{n_{p+1}}. 
\end{equation} 
In fact, 
\begin{align*}
\sum_{p=1}^{n-1} &(\Sd a_p) Q^{(p)}_{n-p}(a) = 
\sum_{p=1}^{n-1} \left[ -a_p - \sum_{k=1}^{p-1} 
\underset{p_1,\dots,p_{k+1}>0}{\sum_{p_1+\cdots+p_{k+1}=p}} (-1)^k
\ll(p_1,\dots,p_k)\ a_{p_1}\cdots a_{p_{k+1}} \right] \\ 
&\kern3cm 
\times 
\left[\sum_{l=1}^{n-p} \underset{j_1,\dots,j_l>0}{\sum_{j_1+\cdots+j_l=n-p}} 
\binom{p+1}{l}\ a_{j_1}\cdots a_{j_l} \right] \\ 
&= \sum_{p=1}^{n-1} \sum_{l=1}^{n-p} 
\underset{j_1,\dots,j_l>0}{\sum_{j_1+\cdots+j_l=n-p}} 
(-1) \binom{p+1}{l}\ a_p a_{j_1}\cdots a_{j_l} \\ 
&\kern.5cm 
+ \sum_{p=1}^{n-1} \sum_{k=1}^{p-1} \sum_{l=1}^{n-p} 
\underset{j_1+\cdots+j_l=n-p}{\sum_{p_1+\cdots+p_{k+1}=p}} 
(-1)^{k+1} \ll(p_1,\dots,p_k) \binom{p+1}{l}
\ a_{p_1}\cdots a_{p_{k+1}} a_{j_1}\cdots a_{j_l} \\ 
&= \sum_{q=1}^{n-1} \sum_{p=1}^{n-q} 
\underset{j_1,\dots,j_q>0}{\sum_{p+j_1+\cdots+j_q=n}} 
(-1) \binom{p+1}{q}\ a_p a_{j_1}\cdots a_{j_q} \\ 
&\kern.5cm 
+ \sum_{q=1}^{n-1} \sum_{k=1}^{q-1} 
\underset{j_1+\cdots+j_{q-k}=n}{\sum_{p_1+\cdots+p_{k+1}}} 
(-1)^{k+1} \ll(p_1,\dots,p_k) \binom{p_1+\cdots +p_{k+1}+1}{q-k}
\ a_{p_1}\cdots a_{p_{k+1}} a_{j_1}\cdots a_{j_{q-k}} \\ 
&= \sum_{q=1}^{n-1}
\underset{n_1,\dots,n_{q+1}>0}{\sum_{n_1+\cdots+n_{q+1}=n}} 
\ a_{n_1}\cdots a_{n_{q+1}} \\ 
&\kern2cm
 \times \left[ -\binom{n_1+1}{q} 
+\sum_{k=1}^{q-1} (-1)^{k+1} \ll(n_1,\dots,n_k) 
\binom{n_1+\cdots +n_{k+1}+1}{q-k} \right] . 
\end{align*} 
Therefore, the identity \eqref{halfantipode} is verified if and only if 
for any $q=1,\dots,n-1$, and for any positive integers $n_1,\dots,n_{q+1}$
with constant sum $n_1+\cdots+n_{q+1}=n$, we have 
\begin{equation*}
(-1)^q \ll(n_1,\dots,n_q) = -\binom{n_1+1}{q} 
+\sum_{k=1}^{q-1} (-1)^{k+1} \ll(n_1,\dots,n_k) 
\binom{n_1+\cdots +n_{k+1}+1}{q-k}. 
\end{equation*} 
This identity is proved in the next lemma.
(Recall the definition \eqref{ll} of the coefficients
$\ll(n_1,\dots,n_q)$.) 
\end{proof} 

\begin{lemma} 
Let $n\geq 2$, then for any $q=1,\dots,n-1$, and for any positive
integers $n_1,\dots,n_{q+1}$,
we have 
\begin{equation}
\label{lemmaeq}
-\binom {n_1+1}{q} + 
\sum _{k=1} ^{q}(-1)^{k+1}\underset{h=1,\dots,k-1}{\underset {
m_1+\dots+m_h\ge h}{\sum_{m_1+\dots+m_k=k}}}
\binom {n_1+1}{m_1}\cdots \binom {n_k+1}{m_k}
\binom {n_1+\dots+n_{k+1}+1}{q-k} = 0.
\end{equation}
\end{lemma}

\begin{proof}
Let us introduce some short notations: let 
$M_j:=m_1+\dots+m_j$, $N_j:=n_1+\dots+n_j$, 
denote by ${\sum}'_{m_1,m_2,\dots}$ the sum over non-negative
integers $m_1,m_2,\dots$ such that $m_1+\dots+m_h\ge h$ for all $h$,
$1\le h\le k-1$, 
and finally set $\Pi_0:=1$ and 
\begin{equation*}
\Pi_j := \prod _{i=1} ^{j}\binom {n_i+1}{m_i}, \quad j>0. 
\end{equation*} 
Using these notations, Equation~\eqref{lemmaeq} may be rewritten in
the form
\begin{equation*}
\sum _{k=0}
^{q}(-1)^{k+1}
{\sum_{M_k=k}\kern-4pt}{\vphantom{\Big\vert}}'\kern4pt
\Pi_{k}
\binom {N_{k+1}+1}{q-k} = 0.
\end{equation*}

Let $S(q)$ be the sum on the left-hand side. We claim that 
for $0\le m\le q$ we have
\begin{equation}
\label{eq}
S(q) = \sum _{k=0} ^{q-m}(-1)^{k+1}
{\sum_{M_k=k}\kern-4pt}{\vphantom{\Big\vert}}'\kern4pt
\Pi_{k} \binom {N_{k+1}+1}{q-k} + \sum _{\ell=1} ^{m}(-1)^{q+\ell}
{\sum _{M_{q+1-m}=q+1-\ell} ^{}\kern-20pt}{\vphantom{\Big\vert}}'\kern10pt
\Pi_{q+1-m} \binom {N_{q+1-m}+\ell-m}{\ell-1}.
\end{equation} 
We prove this claim by induction on $m$. Clearly, formula \eqref{eq} holds
for $m=0$. So, let us suppose that it holds for $m$, and let us, under
this hypothesis, do the following computation:
\begin{align*} 
S(q)&= 
\sum _{k=0}
^{q-m-1}(-1)^{k+1}
{\sum_{M_k=k}\kern-4pt}{\vphantom{\Big\vert}}'\kern4pt
\Pi_{k}
\binom {N_{k+1}+1}{q-k}\\
&\kern1cm+
(-1)^{q-m+1}
{\sum_{M_{q-m}=q-m}\kern-15pt}{\vphantom{\Big\vert}}'\kern10pt
\Pi_{q-m}
\binom {N_{q-m+1}+1}{m}\\
&\kern1cm+
\sum _{\ell=1} ^{m}(-1)^{q+\ell}
\sum _{s=0} ^{m-\ell+1}
{\sum _{M_{q+1-m}=q+1-\ell-s} ^{}\kern-20pt}{\vphantom{\Big\vert}}'\kern10pt
\Pi_{q-m}\binom {n_{q-m+1}+1}s
\binom {N_{q+1-m}+\ell-m}{\ell-1}\\
&=
\sum _{k=0}
^{q-m-1}(-1)^{k+1}
{\sum_{M_k=k}\kern-4pt}{\vphantom{\Big\vert}}'\kern4pt
\Pi_{k}
\binom {N_{k+1}+1}{q-k}\\
&\kern1cm+
\sum _{\ell=1} ^{m+1}(-1)^{q+\ell}
\sum _{s=0} ^{m-\ell+1}
{\sum _{M_{q-m}=q+1-\ell-s} ^{}\kern-20pt}{\vphantom{\Big\vert}}'\kern10pt
\Pi_{q-m}\binom {n_{q-m+1}+1}s
\binom {N_{q+1-m}+\ell-m}{\ell-1}.
\end{align*}
We now replace $\ell$ by $r-s$ and we interchange the inner sums: 
\begin{align} \label{S(q)} 
S(q)&=
\sum _{k=0}
^{q-m-1}(-1)^{k+1}
{\sum_{M_k=k}\kern-4pt}{\vphantom{\Big\vert}}'\kern4pt
\Pi_{k}
\binom {N_{k+1}+1}{q-k}\\
\notag
&\kern0.7cm+
\sum _{r=1} ^{m+1}
{\sum _{M_{q-m}=q+1-r} ^{}\kern-15pt}{\vphantom{\Big\vert}}'\kern10pt
\Pi_{q-m}
\sum _{s=0} ^{r-1}(-1)^{q+r-s}
\binom {n_{q-m+1}+1}s
\binom {N_{q+1-m}+r-s-m}{r-s-1}.
\end{align}
In the second line, the inner sum over $s$ can be evaluated using the
Chu--Vandermonde formula (see e.g.\ \cite[Sec.~5.1, (5.27)]{Graham}). 
Thus, we obtain: 
\begin{align*} 
\sum _{s=0} ^{r-1}(-1)^{q+r-s}&
\binom {n_{q-m+1}+1}s
\binom {N_{q+1-m}+r-s-m}{r-s-1}\\
&=
\sum _{s=0} ^{r-1}(-1)^{q+1}
\binom {n_{q-m+1}+1}s
\binom {-N_{q+1-m}+m-2}{r-s-1}\\
&=
(-1)^{q+1}\binom {n_{q-m+1}-N_{q+1-m}+m-1}{r-1}\\
&=
(-1)^{q+1}\binom {-N_{q-m}+m-1}{r-1}=
(-1)^{q+r}\binom {N_{q-m}+r-m-1}{r-1}.
\end{align*}
If we substitute this in \eqref{S(q)}, 
we obtain exactly formula \eqref{eq} with $m$ replaced by $m+1$.

To prove the lemma, we set $m=q$ in \eqref{eq}. 
This gives 
\begin{align*}
S(q)&=-\binom {n_{1}+1}{q}
+
\sum _{\ell=1} ^{q}(-1)^{q+\ell}
\binom {n_1+1}{q+1-\ell}
\binom {n_{1}+\ell-q}{\ell-1}\\
&= 
\sum _{\ell=1} ^{q+1}(-1)^{q+\ell}
\binom {n_1+1}{q+1-\ell}
\binom {n_{1}+\ell-q}{\ell-1}.
\end{align*}
Again, the sum can be evaluated using the 
Chu--Vandermonde formula. As a result, we obtain
\begin{equation*}
S(q) = (-1)^{q+1}\binom {q-1}q = 0.
\end{equation*}
\end{proof}


\subsection{A tree labelling for the antipode coefficients} 
\label{tree-labelling}

The coefficients $\lambda(n_1,n_2,\dots,n_k)$ in formula
\eqref{antipode} for the antipode are given, by means of \eqref{ll},
as a sum over $k$-tuples $(m_1,m_2,\dots,m_k)$ of non-negative
integers satisfying the two conditions
\begin{align} \label{eq:prop1}
m_1+\dots+m_h&\ge h\quad \text{for all }h=1,2,\dots,k-1,\quad \text{and}\\
\label{eq:prop2}
m_1+\dots+m_k&=k.
\end{align}
Let us denote this set of $k$-tuples by $\mathcal M_k$.
As we are going to outline in this section, it is well known that 
the cardinality of $\mathcal M_k$ is given by the Catalan numbers
$\frac {1} {k+1}\binom {2k}k$. (We refer
the reader to Exercise~6.19 in \cite{Stanley} for 66 combinatorial
interpretations of the Catalan numbers\footnote{with some more recent ones
appearing on Richard Stanley's WWW site
{\tt http://www-math.mit.edu/\~{}rstan/}}, out of which our $k$-tuples
are item~{\bf w.}, modulo the substitutions $n=k$ and $a_i=m_i-1$.)
Thus, in particular, these $k$-tuples are
in bijection with planar binary trees (item~{\bf d.}\ in
\cite[Ex.~6.19]{Stanley}). For the convenience of the
reader, we explain this bijection here in detail.
(What we do, is, essentially, extract the appropriate
restriction of the bijection in \cite[Example~5.3.8]{Stanley}, using
the tree language of Loday \cite{Loday}.)

Recall, from \cite[Sec.~1.5]{Loday} or \cite{BFqedtree}, that for any planar 
binary trees $t, s$, the tree $t$ {\em over} $s$ is defined as the grafting 
$$t \over s := \begin{matrix} \lgraft{$s$}{$t$}\end{matrix}$$ of the root 
of $t$ on the left-most leaf of $s$, and similarly the tree 
$t$ {\em under} $s$ is defined as the grafting 
$$t \under s := \begin{matrix} \rgraft{$t$}{$s$}\end{matrix}$$ 
of the root of $s$ on the right-most leaf of $t$. 
The operations {\em over} and {\em under} are two associative 
(non-commutative) operations with unit given by the ``root tree" $\|$. 
Moreover, any planar binary tree can be written as a monomial in $\Y$, 
with respect to {\em over} and {\em under}. 

Using this notation, the mapping from $\mathcal M_k$ to the set $Y_k$
of planar binary trees with $k$ internal vertices is given by the following 
algorithm.

\begin{definition} \label{def1}
For any $k\geq 1$, let $\Phi: \M_k \longrightarrow Y_k$ be the map defined 
by the following recursive algorithm. 
For any $m=(m_1,\dots,m_k) \in \M_k$, 

\begin{enumerate}
\item if $m=(1)$, then set $\Phi(m):= \Y$; 
\item if $m=(m_1,\dots,m_l,m_{l+1},\dots,m_k) \in \M_k$ is such that 
$(m_1,\dots,m_l) \in \M_l$ and $(m_{l+1},\dots,m_k) \in \M_{k-l}$, then set 
$$\Phi(m):= \Phi(m_1,\dots,m_l) \under \Phi(m_{l+1},\dots,m_k);$$ 
\item if $m=(m_1,\dots,m_{k-1},0) \in \M_k$ is not decomposable in 
sub-tuples as in 2., then $m_1>1$; in this case set
$$\Phi(m):= \Phi(m_1-1,\dots,m_{k-1}) \over \Y.$$ 
\end{enumerate}
\end{definition}

It is easy to show that if a $k$-tuple $m=(m_1,\dots,m_k)$ is not 
decomposable in sub-tuples as in 2., then it is of the form given in 3.
In fact, by Eq.~\eqref{eq:prop2} we can write 
$m_1+\dots+m_{k-1}=k-m_k$ and from Eq.~\eqref{eq:prop1} we get 
$m_k \leq 1$. But $m_k =1$ implies that $(m_k,\dots,m_{k-1}) \in \M_{k-l}$, 
hence $m$ would be decomposable as in 2, a contradiction. 
Therefore $m_k=0$, and 
$m=(m_1,\dots,m_{k_1},0)$. 
Furthermore, if $m_1=1$ the $k$-tuple is decomposable into 
$(m_1) = (1) \in \M_1$ and $(m_2,\dots,m_{k-1},0) \in \M_{k-1}$, which 
contradicts again our original assumption.
We must therefore have $m_1>1$ in this case. 

To give an example, consider the sequence $m=(4,0,1,0,0,2,1,0) \in \M_8$. 
This $8$-tuple is decomposable into the two indecomposable tuples 
$(4,0,1,0,0) \in \M_5$ and $(2,1,0) \in \M_3$, so 
$\Phi(4,0,1,0,0,2,1,0) = \Phi(4,0,1,0,0) \under \Phi(2,1,0)$. 
We now apply Step~3. of the algorithm in Definition~\ref{def1} to the two
terms on the right-hand side separately:
\begin{align*}
\Phi(4,0,1,0,0) & =\Phi(3,0,1,0) \over \Y 
= \Phi(2,0,1) \over \Y \over \Y \\ 
& = \left(\Phi(2,0) \under \Phi(1)\right) \over \Y \over \Y 
= \left(\left(\Phi(1) \over \Y \right) \under \Y \right) \over \Y \over \Y \\ 
& = \left(\left(\Y \over \Y \right) \under \Y \right) \over \Y \over \Y , 
\end{align*}
and 
\begin{align*}
\Phi(2,1,0) & =\Phi(1,1) \over \Y 
= \big(\Phi(1) \under \Phi(1)\big) \over \Y 
=  \Big(\Y \under \Y\Big) \over \Y .
\end{align*}
In conclusion, 
\begin{align*}
\Phi(4,0,1,0,0,2,1,0) 
& = \left[\left(\left(\Y \over \Y \right) \under \Y \right) 
\over \Y \over \Y\right] \under 
\left[ \Big(\Y \under \Y\Big) \over \Y \right] \\  
&\raise1cm\hbox{${}={}$} 
\hbox{\hskip4cm}
\Einheit.5cm
\unitlength\Einheit
\thicklines
\raise0 \Einheit\hbox to0pt{\hskip0 \Einheit\line(-2,1){2}\hss}
\raise0 \Einheit\hbox to0pt{\hskip0 \Einheit\line(2,1){2}\hss}
\raise1 \Einheit\hbox to0pt{\hskip-2 \Einheit\line(-2,1){2}\hss}
\raise1 \Einheit\hbox to0pt{\hskip-2 \Einheit\line(1,1){1}\hss}
\raise1 \Einheit\hbox to0pt{\hskip2 \Einheit\line(-1,1){1}\hss}
\raise1 \Einheit\hbox to0pt{\hskip2 \Einheit\line(1,1){1}\hss}
\raise2 \Einheit\hbox to0pt{\hskip-4 \Einheit\line(-2,1){2}\hss}
\raise2 \Einheit\hbox to0pt{\hskip-4 \Einheit\line(1,1){1}\hss}
\raise2 \Einheit\hbox to0pt{\hskip1 \Einheit\line(-1,1){1}\hss}
\raise2 \Einheit\hbox to0pt{\hskip1 \Einheit\line(1,1){1}\hss}
\raise3 \Einheit\hbox to0pt{\hskip-6 \Einheit\line(-1,1){1}\hss}
\raise3 \Einheit\hbox to0pt{\hskip-6 \Einheit\line(1,1){1}\hss}
\raise3 \Einheit\hbox to0pt{\hskip-3 \Einheit\line(-1,1){1}\hss}
\raise3 \Einheit\hbox to0pt{\hskip-3 \Einheit\line(1,1){1}\hss}
\raise3 \Einheit\hbox to0pt{\hskip2 \Einheit\line(-1,1){1}\hss}
\raise3 \Einheit\hbox to0pt{\hskip2 \Einheit\line(1,1){1}\hss}
\end{align*}

The inverse mapping is given by the algorithm described below.

\begin{definition}
For any $k\geq 1$, let $\Psi: Y_k \longrightarrow \M_k$ be the map defined 
by the following recursive algorithm. For any $t \in Y_k$, 
\begin{enumerate}
\item if $t=\Y$, then set $\Psi(\Y):= (1)$; 
\item if $t= t_1 \under t_2$, then set 
$\Psi(t):= \left(\Psi(t_1),\Psi(t_2)\right)$; 
\item if $t=t_1 \over \Y$, then set 
$\Psi(t):= (m_1+1,\dots,m_k,0)$, where $(m_1,\dots,m_k) = \Psi(t_1)$. 
\end{enumerate}
\end{definition}

\noindent
It is easy to see that always exactly one of the cases 1., 2., or
3.\ applies.


\section{Co-action and semi-direct co-product of the Hopf algebras}
\label{semi-direct}

Since the group $\Gd$ of formal diffeomorphisms acts by composition 
on the group $\Gi$ of invertible series, one can consider the semi-direct 
product $\Gd \ltimes \Gi$ of the two groups. 
If we consider series with non-commutative coefficients, of course the 
semi-direct product $\Gd \ltimes \Gi(A)$ is still a group, while the 
semi-direct product $\Gd(A) \ltimes \Gi(A)$ is not anymore a group, because 
$\Gd(A)$ itself is not a group. 

In this section, we show that the dual construction on the co-ordinate rings 
still makes sense on the non-commutative algebras. It produces a Hopf algebra 
$\C(\Gd) \ltimes \Hi$ which is neither commutative nor co-commutative
in the case corresponding to the semi-direct product group, 
and an algebra $\Hd \ltimes \Hi$ which is also a co-algebra but not 
a bi-algebra in the more general case. 

\subsection{Action and semi-direct product of the groups of series} 

The composition $f \circ g$ of two invertible (formal power) 
series is not a formal power series,  
because the constant term 
$(f \circ g)_0 = \sum_{n=0}^\infty f_n (g_0)^n$ 
is an infinite sum. 
However, an invertible series $f$ can be composed with a formal 
diffeomorphisms $\varphi$, and the result $f \circ \varphi$ is again 
an invertible series. Moreover, we have 
$(f \circ \varphi)\circ \psi = f \circ (\varphi \circ \psi)$.
In other words, the composition $\circ: \Gi \times \Gd \longrightarrow \Gi$ 
is a natural right action of $\Gd$ on $\Gi$. 
Furthermore, this action commutes with the group structure of $\Gi$, 
in the sense that 
\begin{equation}
\label{Gi-action} 
        (f g) \circ \varphi = (f \circ \varphi)\  (g \circ \varphi). 
\end{equation}

In such a situation, we can consider the semi-direct product 
$\Gd \ltimes \Gi$, which is the group defined on the direct product 
$\Gd \times \Gi$ by the law 
\begin{equation}
\label{semidirect product} 
(\varphi,f) \ldot (\psi,g) 
:= \big(\varphi \circ \psi,(f \circ \varphi) g \big). 
\end{equation}

In the dual context, on the co-ordinate rings, we have the co-action 
$\ddif: \C(\Gi) \longrightarrow \C(\Gi) \otimes \C(\Gd)$ which satisfies
$\<b_n,f \circ \varphi \> = \<  \ddif b_n ,f \otimes \varphi\>$,
where the $b_n$'s are the generators of $\C(\Gi)$ defined in
Section~\ref{invertible}. 
As we did in Section~\ref{diffeomorphisms} for the composition 
of formal diffeomorphisms, we can compactly encode the co-action $\ddif$ 
on the generators $b_n$ by introducing the generating series 
\begin{equation}
\label{B(x)} 
B(x) = 1 + \sum_{n=1}^\infty b_n x^n =
\sum_{n=0}^\infty b_n x^n \quad (b_0 :=1) , 
\end{equation}
so that $\varphi(x)=\<B(x),\varphi\>$.
As for the co-product of $\C(\Gd)$, the co-action is then given by the 
formal residue 
\begin{equation}
\label{B(x)-coaction} 
\ddif B(x) = \res B(z) \otimes \frac{1}{z - A(x)}, 
\end{equation}
where $A(x)$ is the generating series of the generators of 
$\C(\Gd)$ as in \eqref{A(x)}. 
The reader should note that we can also describe the co-product $\Di$ 
dual to the product \eqref{product} directly on the generating series. 
It has the simple expression 
\begin{equation}
\label{B(x)-coproduct}
\Di B(x) = B(x) \otimes B(x).
\end{equation}

The co-action $\ddif$ allows us to construct the co-product on the co-ordinate 
ring $\C(\Gd \ltimes \Gi)$ of the semi-direct product group. 
As an algebra, we have $\C(\Gd \ltimes \Gi) \cong \C(\Gd) \otimes \C(\Gi)$, 
and the co-product dual to the product \eqref{semidirect product} is 
\begin{equation}
\label{semidirect coproduct}
\Dl (a \otimes b) = \Dd(a)\ \big[ (\ddif \otimes \Id) \Di(b) \big], 
\qquad a \in \C(\Gd),\ b \in \C(\Gi). 
\end{equation}


\subsection{Dual non-commutative co-action and semi-direct co-product}
\label{co-action}

The previous discussion can be repeated for the group $\Gi(\A)$ of 
invertible series with coefficients in $\A$, and for its dual 
non-commutative Hopf algebra $\Hi = \C\<b_1,b_2,\dots\>$. 
If we adopt the definition \eqref{B(x)} of the non-commutative 
generating series $B(x)$, and if $A(x)$ denotes the 
non-commutative generating series 
for the generators of the 
Hopf algebra $\Hd$, then we can still use formula \eqref{B(x)-coaction} 
to define a co-action $\ddif: \Hi \longrightarrow \Hi \otimes \Hd$ 
of $\Hd$ on $\Hi$. 

\begin{lemma} \label{delta}
The explicit expression of the co-action $\ddif$ on the generators of $\Hi$ 
is: \begin{align*}
\ddif b_n &=  \sum_{k=0}^{n} b_k \otimes Q^{(k-1)}_{n-k}(a),\quad n\ge0,
\end{align*}
where we use the identification $b_0=1$, and where the polynomials
$Q^{(k)}_m(a)$ are the polynomials from Definition~{\em\ref{defQ}}.
\end{lemma}

\begin{proof} 
We compute the explicit expression for the co-action by applying it
to $B(x)$:
\begin{align}
\ddif B(x) &= \res \sum_{n=0}^\infty \sum_{k=0}^\infty 
b_k \otimes A(x)^n z^{k-n-1} 
= \sum_{n=0}^\infty b_n \otimes A(x)^n \label{dddB(x)}\\ 
&= \sum_{n=0}^\infty x^n \bigg( 
\sum_{k=0}^{n} b_k \otimes Q^{(k-1)}_{n-k}(a) \bigg). \notag
\end{align}
This proves the formula in the statement of the lemma. 
\end{proof}

The first few values of the co-action on the generators $b_i$ are
\begin{align*}
\ddif b_1 &= b_1 \otimes 1, \\
\ddif b_2 &= b_2 \otimes 1 + b_1 \otimes a_1 , \\
\ddif b_3 &= b_3 \otimes 1 + 2 b_2 \otimes a_1 + b_1 \otimes a_2 .
\end{align*}

\begin{lemma} 
The map $\ddif: \Hi \longrightarrow \Hi \otimes \Hd$ is a co-action with 
respect to $\Dd$, that is  
\begin{equation*}
(\ddif \otimes \Id) \ddif = (\Id \otimes \Dd) \ddif. 
\end{equation*}
\end{lemma} 

\begin{proof} 
It suffices to prove this equality for the generators $b_n$, or,
equivalently, 
it suffices to prove it for the generating series $B(x)$:  
the left-hand side, by formula \eqref{B(x)-coaction}, is
\begin{align*}
(\ddif \otimes \Id) \ddif B(x) 
&= \resz{z_1} \ddif B(z_1) \otimes \frac{1}{z_1 - A(x)} \\
&= \resz{z_1} \resz{z_2} B(z_2) \otimes \frac{1}{z_2 - A(z_1)} 
\otimes \frac{1}{z_1 - A(x)}, 
\end{align*}
and the right-hand side, by formulae \eqref{B(x)-coaction} and 
\eqref{Ddfrac}, is
\begin{align*}
(\Id \otimes \Dd) \ddif B(x) &= \resz{z_2} 
B(z_2) \otimes \Dd \big(\frac{1}{z_2 - A(x)}\big) \\
&= \resz{z_2} \resz{z_1} B(z_2) \otimes \frac{1}{z_2 - A(z_1)} 
\otimes \frac{1}{z_1 - A(x)}.
\end{align*}
Thus, the co-action is co-associative. 
\end{proof}  

\begin{lemma} 
The Hopf algebra $\Hi$ is a co-algebra comodule over $\Hd$ (in the sense of 
\cite{Molnar}), that is 
\begin{equation*}
(\Di \otimes \Id)\ \ddif = (\Id \otimes \Id \otimes m) 
(\Id \otimes \tau \otimes \Id)(\ddif \otimes \ddif)\ \Di, 
\end{equation*} 
where $\tau$ is the twist operator $\tau(u \otimes v) = v \otimes u$ and 
$m$ denotes the multiplication in $\Hd$. 
\end{lemma}

\begin{proof} 
We check this identity on the generating series $B(x)$.
The co-product $\Di$ on the generating series is given by 
\eqref{B(x)-coproduct}, and by \eqref{Gp-coproduct} on each $b_i$, 
while the co-action is given by \eqref{dddB(x)}. 
Thus, we have 
\begin{align*}
(\Di\otimes \Id)\ \ddif B(x) &= 
(\Di\otimes \Id)\ \sum_{k=0}^\infty b_k \otimes A(x)^k \\ 
&= \sum_{k=0}^\infty \sum_{m+n=k} b_m \otimes b_n \otimes A(x)^k    
= \sum_{m,n=0}^\infty b_m \otimes b_n \otimes A(x)^{m+n} \\ 
&= (\Id \otimes \Id \otimes m) (\Id \otimes \tau \otimes \Id) 
\big( \sum_{m=0}^\infty b_m \otimes A(x)^m \otimes 
\sum_{n=0}^\infty b_n \otimes A(x)^n \big) \\ 
&= (\Id \otimes \Id \otimes m) (\Id \otimes \tau \otimes \Id) 
(\ddif B(x) \otimes \ddif B(x)) \\ 
&= (\Id \otimes \Id \otimes m) (\Id \otimes \tau \otimes \Id) 
(\ddif \otimes \ddif) \Di B(x). 
\end{align*}
\end{proof} 

By Molnar's results \cite{Molnar}, in such a situation we can 
consider the semi-direct or smash co-product $\Hd \ltimes \Hi$ 
of Hopf algebras. This space is at the same time an algebra and a co-algebra, 
with co-product given by formula \eqref{semidirect coproduct}. 
However, the co-product is not an algebra homomorphism, because $\Hd$ 
is not a commutative algebra, and therefore $\Hd \ltimes \Hi$ is not 
a bi-algebra. In order to obtain a Hopf algebra structure, we should 
consider the semi-direct co-product $\C(\Gd) \ltimes \Hi$ constructed 
in the same way. Despite the commutativity of $\C(\Gd)$ and the 
co-commutativity of $\Hi$, the resulting Hopf algebra is neither 
commutative nor co-commutative. 


\section{Relation with the QED renormalization Hopf algebras} 
\label{QEDHopf}

In \cite{BFqedtree} and \cite{BFqedren}, it was shown that the 
renormalization of the electron propagator in quantum electrodynamics 
can be described in terms of a semi-direct co-product Hopf algebra 
$\Ha \ltimes \He$ on the set of rooted planar binary trees. 
Here, $\Ha$ is a commutative Hopf algebra which represents the renormalization 
of the electric charge, while $\He$ is a Hopf algebra which represents 
the electron propagators, and which is neither commutative nor co-commutative. 
The renormalization is then a co-action of $\Ha \ltimes \He$ on $\He$, 
obtained as a restriction of the co-product.  

The Hopf algebra $\Ha$ also describes the renormalization of the photon 
propagators, by means of a co-action of $\Ha$ on the non-commutative Hopf 
algebra $\Hp$ which represents the photon propagators. It turns out that the 
generators of $\Hp$ are exactly all the elements of $\Ha$, and that 
the co-action on the generators of $\Hp$ coincides with the co-product 
in $\Ha$. 
Since $\Hp$ is a non-commutative algebra, a suitable restriction of 
this co-action can be interpreted as a co-product on a non-commutative 
extension $\Hat$ of the commutative Hopf algebra $\Ha$. 

In this section, we show that $\Hi$ is a Hopf sub-algebra of both $\He$ 
and $\Hp$, and that $\Hd$ is a Hopf sub-algebra of $\Hat$. 
By construction, it follows also that the co-ordinate ring $\C(\Gd)$ 
is a Hopf sub-algebra of $\Ha$, and that the semi-direct co-product 
$\C(\Gd) \ltimes \Hi$ of Section~\ref{semi-direct} is a Hopf sub-algebra of 
the QED Hopf algebra $\Ha \ltimes \He$. 


\subsection{$\Hd$ and the charge renormalization Hopf algebra $\Hat$ on trees} 

We start by recalling the definition of $\Hat$ from \cite{BFqedtree}.
As in Section~\ref{tree-labelling}, denote by $Y_n$ the set of rooted planar 
binary trees with $n$ internal vertices, and by $Y_\infty=\bigcup_{n\geq0}Y_n$ 
the union of all such trees. 
Finally, denote by $\Hat:=\C Y_\infty$ the vector space spanned by all trees. 
Then $\Hat$ is a non-commutative unital algebra with the {\em over} product 
$t \over s$ of Section~\ref{tree-labelling}, and unit given by $\|$. 
In particular, $\Hat$ is a graded and connected algebra, with 
graded components $(\Hat)_n = \C Y_n$. 

Moreover, $\Hat$ is a free algebra (i.e., a tensor algebra). 
To see this, we denote 
by $\vee: Y_n \times Y_m \longrightarrow Y_{n+m+1}$ the map which grafts 
two trees on a new root. If, for any $t \in Y_n$, we set 
$V(t)= \| \vee t \in Y_{n+1}$, then $\Hat $ 
is isomorphic to the free algebra
$ \C\<V(t),t\in Y_\infty\>$ generated by the 
trees $V(t)$. Indeed, any tree $t$ can be 
decomposed as 
$$t = t^l \vee t^r = t^l \over V(t^r) = \cdots = 
V(t^{ll\dots l})\over V(t^{ll\dots r})
\over\cdots\over V(t^{lr})\over V(t^r).$$ 
As usual, we identify $\|$ with the element $1$ in
$ \C\<V(t),t\in Y_\infty\>$. 

In \cite{BFqedtree}, it was
shown that $\Hat$ is a connected Hopf algebra
which is neither commutative nor co-commu\-ta\-tive. 
We recall briefly the explicit definitions. 
The co-product $\Da: \Hat \longrightarrow \Hat \otimes \Hat$ is defined 
recursively by the formulae 
\begin{eqnarray*}
\Da \| &=& \| \otimes \|, \\ 
\Da V(r) &=& \| \otimes V(r) + \da V(r), \\ 
\Da (r \vee s) &=& \Da r \over \Da V(s) ; 
\end{eqnarray*}
where $\da : \Hat \longrightarrow \Hat \otimes \Hat$ is the right 
co-action\footnote{
A right co-action $\da$ of $\Hat$ on itself satisfies the co-associativity 
condition $(\da \otimes \Id ) \da = (\Id \otimes \Da) \da$.} 
of $\Hat$ on itself given by the recursive formulae 
\begin{eqnarray*}
\da \| &=& \| \otimes \|, \\  
\da V(r) &=& (V \otimes \Id) \da (r), \\  
\da (r \vee s) &=& \Da r \over \da (V(s)) .  
\end{eqnarray*}
The co-unit $\varepsilon : \Hat \longrightarrow \C$ is the linear map 
which sends all the trees to $0$, except for the ``root tree'' $\|$ 
which is sent to $1$. 
The antipode is defined by a standard recursive formula similar to
the one in \eqref{recform}, since the algebra $\Hat$ is connected. 

In the statement of the theorem below, we need one more notation: we
write $|t|$ for the number of internal vertices of the tree $t$. 

\begin{theorem} 
\label{HdHat}
The map $\Hd \longrightarrow \Hat$ given by 
\begin{align*}
\Omega:a_n \mapsto t_n:=\sum_{|t|=n}t, 
\end{align*} 
and extended as a homomorphism of unital algebras,
is an injective co-algebra homomorphism. In particular, the Hopf algebra $\Hd$ 
is a Hopf sub-algebra of $\Hat$. 
\end{theorem}

\begin{proof} 
To prove that $\Omega$ is injective, we consider a
(non-commutative) polynomial $P(a)$ in the $a_n$'s in the kernel of
$\Omega$, that is, we have $P(t)=0$, where $P(t)$ is obtained from
$P(a)$ by replacing $a_n$ by $t_n$ for all $n$ and the product by the
over product $\over$. Obviously, $P(a)$ has no
constant term, because otherwise we would trivially have $P(t)\ne0$. 
Let us suppose that $P(a)$ is not identically zero. 
In that case, there exists a monomial
$a_{i_1}a_{i_2}\cdots a_{i_k}$ which appears with non-zero coefficient
in $P(a)$. Without loss of generality, we may assume that $k$ is
minimal. In particular, since $P(a)$ has no constant term, we have
$k\ge1$. 

\begin{figure}
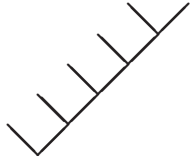

$$
\Einheit.4cm
\Pfad(0,0),33333\endPfad
\Pfad(1,1),8\endPfad
\Pfad(2,2),8\endPfad
\Pfad(3,3),8\endPfad
\Pfad(4,4),8\endPfad
\Pfad(0,0),8\endPfad
\hskip1.5cm
$$
\caption{The right brush $r_5$}
\label{fig:r}
\end{figure}

\begin{figure}
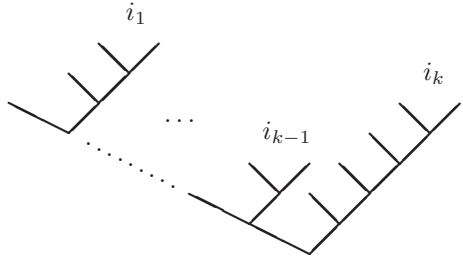

$$
\Einheit.4cm
\Pfad(0,0),33333\endPfad
\Pfad(1,1),8\endPfad
\Pfad(2,2),8\endPfad
\Pfad(3,3),8\endPfad
\Pfad(4,4),8\endPfad
\unitlength\Einheit
\thicklines
\raise0 \Einheit\hbox to0pt{\hskip0 \Einheit\line(-2,1){2}\hss}
\Pfad(-2,1),33\endPfad
\Pfad(-1,2),8\endPfad
\thicklines
\raise1 \Einheit\hbox to0pt{\hskip-2 \Einheit\line(-2,1){2}\hss}
\Pfad(-8,4),333\endPfad
\Pfad(-6,6),8\endPfad
\Pfad(-7,5),8\endPfad
\thicklines
\raise4 \Einheit\hbox to0pt{\hskip-8 \Einheit\line(-2,1){2}\hss}
\raise25pt\hbox to0pt{\hskip-53pt.\hss}
\raise27pt\hbox to0pt{\hskip-57pt.\hss}
\raise29pt\hbox to0pt{\hskip-61pt.\hss}
\raise31pt\hbox to0pt{\hskip-65pt.\hss}
\raise33pt\hbox to0pt{\hskip-69pt.\hss}
\raise35pt\hbox to0pt{\hskip-73pt.\hss}
\raise37pt\hbox to0pt{\hskip-77pt.\hss}
\raise39pt\hbox to0pt{\hskip-81pt.\hss}
\raise41pt\hbox to0pt{\hskip-85pt.\hss}
\raise50pt\hbox to0pt{\hskip-55pt\dots\hss}
\Label\l{i_1}(-5,8)
\Label\l{i_{k-1}}(-0,4)
\Label\r{\hbox{\kern-14pt$i_{k}$}}(4,6)
\hskip-2cm
$$
\caption{The over product $r_{i_1}\over r_{i_2}\over \cdots\over r_{i_k}$}
\label{fig:rrr}
\end{figure}

Le us now consider the image of this monomial under $\Omega$,
\begin{equation} \label{eq:tmon} 
t_{i_1}\over t_{i_2}\over \cdots\over t_{i_k}.
\end{equation}
For any $i$, in the expansion of
$t_{i}$, there appears the ``right brush" $r_{i}$, which, by
definition, is the planar binary tree consisting of $i$ internal
nodes, all of which (except for the root) are right descendants of
another internal node. (See Figure~\ref{fig:r} for an illustration of
the right brush $r_5$).
Hence, in the expansion of the monomial \eqref{eq:tmon}, there appears
the over product $r_{i_1}\over r_{i_2}\over \cdots\over r_{i_k}$ (see
Figure~\ref{fig:rrr}). 
As is not difficult to see, this tree cannot appear in the
expansion of any other monomial $t_{j_1}\over t_{j_2}\over \cdots
\over t_{j_l}$ with
$l\ge k$. Therefore, it cannot cancel out in the expansion of $P(t)$, a
contradiction.

\medskip
To prove that $\Hd$ is a sub-co-algebra of $\Hat$, we only need to show 
that $\Da t_n=(\Omega\otimes\Omega)\Dd a_n$. To do this, we actually 
prove also that $\Hd$ is a right sub-co-module of $\Hat$, where the 
right co-action of $\Hd$ on itself is induced by the natural right 
co-action $\ddif: \Hi \longrightarrow \Hi \otimes \Hd$ defined in 
Section~\ref{co-action}, via the co-module homomorphism 
$\Hd \longrightarrow \Hi$, $a_n \mapsto b_n$. 
To be precise, we are going to show that 
\begin{align} 
\Da(t_n)&= 
 \sum_{k=0}^{n} t_{k} \otimes Q^{(k)}_{n-k}(t), \label{Da(an)}\\ 
\da(t_n)&= 
 \sum_{k=0}^{n} t_{k} \otimes Q^{(k-1)}_{n-k}(t), \label{da(an)}
\end{align}
where the polynomial $Q^{(m)}_{n}(t)$ is as in 
Definition~\ref{defQ}, only that the $a_n$'s get replaced by the $t_n$'s
and the product by the over product $\over$. 
By Lemmas~\ref{Dd(a_n)} and \ref{delta}, this would accomplish the
proof that $\Hd$ is a sub-co-algebra of $\Hat$.

We prove these two claims simultaneously
by induction on $n$. 
To simplify notation in the following calculations, we shall from now
on omit the symbol $/$ in the {\it over} product $x\over y$ of the elements 
$x$ and $y$ of $\Hat$ and write simply $xy$.
 
For $n=0$ we have 
\begin{align*} 
\Da(\|) &= \| \otimes \|, \\
\da(\|) &= \| \otimes \|, 
\end{align*}
while for $n=1$ we have 
\begin{align*} 
\Da(\Y)&=\Y\otimes\|+\|\otimes\Y,\\
\da(\Y) &= \Y \otimes \|.
\end{align*}
So, formulae \eqref{Da(an)} and \eqref{da(an)} hold for $n=0$ and
$n=1$ because 
$t_0 = \|$ and $t_1 = \Y$. 
Now we suppose that they hold up to a fixed $n \geq 1$, and we show that 
they hold for $n+1$. 

For this induction step, we repeatedly need an expansion formula for
$t_{n+1}$. 
If we decompose $t=r \vee s$ into a left tree $r$ and a right tree
$s$, with $|t|=|r|+|s|+1$, then we have $t=r\over V(s)=rV(s)$, and thus
\begin{align} 
t_{n+1} &= \underset{|r|+|s|=n}{\sum_{0 \leq |r|,|s| \leq n}} r  V(s) 
 \nonumber \\ 
&= \sum_{m=0}^{n} t_{n-m}  V(t_m)  
\label{a{n+1}}. 
\end{align}

\medskip
{\it Induction step for \eqref{Da(an)}.} 
Using the definitions of $\Da$ and $\da$,
for any $m \leq n$ we have
\begin{align} 
\Da V(t_m) &= \sum_{|t|=m} \Da V(t) 
= \sum_{|t|=m} \|\otimes V(t) + \sum_{|t|=m} \da V(t) \nonumber \\ 
&= \sum_{|t|=m} \|\otimes V(t) + \sum_{|t|=m} (V\otimes\Id) \da(t) \nonumber\\ 
&= \| \otimes V(t_m) 
+ \sum_{k=0}^{m} V(t_{k}) \otimes Q_{m-k}^{(k-1)}(t).
\label{Da(Vak)} 
\end{align}
Therefore, using Eq.~\eqref{a{n+1}}, the definition of $\Da$, and 
Eqs.~\eqref{Da(an)} and \eqref{Da(Vak)}, we get 
\begin{align} 
\notag
\Da(t_{n+1}) &= 
\sum_{m=0}^n \Da(t_{n-m})  \Da V(t_m) \\ 
&= 
 \sum_{m=0}^{n} \sum_{k=0}^{n-m} 
t_k \otimes Q_{n-m-k}^{(k)}(t)  V(t_{m})
+ \sum_{m=0}^{n} \sum_{k=0}^{n-m} \sum_{l=0}^{m} 
t_k  V(t_l)\otimes Q_{n-m-k}^{(k)}(t)  Q_{m-l}^{(l-1)}(t).
\label{twoterms}
\end{align}
We simplify the two terms on the right-hand side separately.
In the first term, we use the recurrence for $ Q_{n-m-k}^{(k)}(t)$ 
given in Lemma~\ref{recursiveQ}, Eq.~\eqref{a{n+1}}, and then
again Lemma~\ref{recursiveQ}, to obtain
\begin{align} \notag
 \sum_{m=0}^{n} \sum_{k=0}^{n-m} 
t_k \otimes Q_{n-m-k}^{(k)}(t)  V(t_{m})
&= \sum_{m=0}^{n} \sum_{k=0}^{n-m} 
t_k \otimes
\sum _{l=0} ^{n-m-k} Q_{n-m-k-l}^{(k-1)}(t)t_l  V(t_{m})\\
\notag
&=  \sum_{k=0}^{n} \sum_{p=0}^{n-k}
t_k \otimes
 Q_{n-k-p}^{(k-1)}(t)\sum _{l=0} ^{p}t_l  V(t_{p-l})
\\
\notag
&=  \sum_{k=0}^{n} \sum_{p=0}^{n-k}
t_k \otimes
 Q_{n-k-p}^{(k-1)}(t)t_{p+1}
\\
\notag
&= \sum_{k=0}^{n}   
t_k \otimes \sum_{i=0}^{n-k+1}
 Q_{n-k-i+1}^{(k-1)}(t)t_{i}
-\sum_{k=0}^{n} 
t_k \otimes
 Q_{n-k+1}^{(k-1)}(t)\\
\label{term1}
&= \sum_{k=0}^{n}   
t_k \otimes 
 Q_{n-k+1}^{(k)}(t)
-\sum_{k=0}^{n} 
t_k \otimes
 Q_{n-k+1}^{(k-1)}(t).
\end{align}
On the other hand, to the second term in \eqref{twoterms}
we apply 
the quadratic identity satisfied by the $Q_n^{(m)}(t)$'s proved in 
Lemma~\ref{quadraticQ} and \eqref{a{n+1}}, to obtain 
\begin{align}
\notag
\sum_{k=0}^{n} \sum_{l=0}^{n-k} 
t_k  V(t_l)\otimes \Bigg(
 \sum_{m=l}^{n-k} Q_{n-m-k}^{(k)}(t)  Q_{m-l}^{(l-1)}(t)\Bigg)
  &=   \sum_{k=0}^{n} \sum_{l=0}^{n-k} 
t_k  V(t_l)\otimes Q_{n-k-l}^{(k+l)}(t) \\
\notag
&= \sum_{p=0}^{n} \sum_{l=0}^{p} 
t_{p-l}  V(t_l)\otimes Q_{n-p}^{(p)}(t) \\
\notag
&= \sum_{p=0}^{n}t_{p+1}
\otimes Q_{n-p}^{(p)}(t) \\
\label{term2}
&= \sum_{p=0}^{n+1}t_{p}
\otimes Q_{n-p+1}^{(p-1)}(t) .
\end{align}
By summing the two expressions \eqref{term1} and \eqref{term2}, we
arrive exactly at \eqref{Da(an)} with $n$ replaced by $n+1$, 
because $Q^{(n)}_0(t)=Q^{(n+1)}_0(t)$.

\medskip
{\it Induction step for \eqref{da(an)}.} 
Using the definition of $\da$, 
for any $m \leq n$, we have
\begin{align} 
\da V(t_m) &= \sum_{|t|=m} \da V(t) 
= \sum_{|t|=m} (V\otimes\Id) \da(t) \nonumber\\ 
&=  \sum_{k=0}^{m} V(t_{k}) \otimes Q_{m-k}^{(k-1)}(t). 
\label{da(Vak)} 
\end{align}
Therefore, using Eq.~\eqref{a{n+1}}, the definition of $\da$,
Eqs.~\eqref{Da(an)} and \eqref{da(Vak)}, we get 
\begin{align*} 
\da(t_{n+1}) &= 
\sum_{ m =0}^n \Da(t_{n-m})  \da V(t_m) \\ 
&= \sum_{m=0}^{n} \sum_{k=0}^{n-m} \sum_{l=0}^{m} 
t_k  V(t_l)\otimes Q_{n-m-k}^{(k)} (t) Q_{m-l}^{(l-1)}(t)\\
&=  \sum_{k=0}^{n} \sum_{l=0}^{n-k} 
t_k  V(t_l)\otimes \Bigg(\sum_{m=l}^{n-k}
Q_{n-m-k}^{(k)} (t) Q_{m-l}^{(l-1)}(t)\Bigg).
\end{align*}
We apply again the quadratic relation of Lemma~\ref{quadraticQ} to the
sum over $m$. This leads to
\begin{align*} 
\da(t_{n+1}) 
&=  \sum_{k=0}^{n} \sum_{l=0}^{n-k} 
t_k  V(t_l)\otimes 
Q_{n-k-l}^{(k+l)} (t) \\
&=  \sum_{p=0}^{n}\Bigg( \sum_{l=0}^{p} 
t_{p-l}  V(t_{l})\Bigg)\otimes 
Q_{n-p}^{(p)} (t). 
\end{align*}
Finally, another application of the recurrence \eqref{a{n+1}} yields
\begin{align*}
\da(t_{n+1}) 
&= \sum_{p=0}^{n}
t_{p+1}\otimes 
Q_{n-p}^{(p)} (t),
\end{align*}
which, after replacing $p$ by $k-1$, agrees with \eqref{da(an)} with
$n$ replaced by $n+1$.
\end{proof}


\subsection{$\Hi$ and the propagator Hopf algebras $\He$ and $\Hp$ on trees} 

Denote by $\He$ and $\Hp$ the free associative algebras 
$\C\<Y_\infty\>/(\|-1)$ on the set of all trees, where we identify the 
``root tree'' $\|$ with the unit. 
As we mentioned in Section~\ref{tree-labelling}, the {\em under} and the 
{\em over} products on trees are associative and have $\|$ as a unit. 
Therefore their dual co-operations are co-unital and co-associative, 
and their multiplicative extensions on $\He$ and $\Hp$, denoted 
respectively by $\DPe$ and $\DPp$, define two structures of a Hopf algebra 
on $\He$, respectively on $\Hp$, 
which are neither commutative nor co-commutative. 

The respective co-products can be 
defined on the generators $t=r\vee s$ in a recursive manner, by putting 
\begin{align}
\label{DPe}
\DPe(r\vee s) &= 1 \otimes (r\vee s) + 
\sum_{} (r\vee s\ix1) \otimes s\ix2,\quad \text{where }
\sum _{} ^{}s\ix1\otimes s\ix2=\DPe(s), \\ 
\label{DPp}
\DPp(r\vee s) &= (r\vee s) \otimes 1 + 
\sum_{} r\ix1 \otimes (r\ix2\vee s), \quad \text{where }
\sum _{} ^{}r\ix1\otimes r\ix2=\DPp(r).
\end{align}
Here we use Sweedler's notation $\sum t\ix1 \otimes t\ix2$ for both 
$\DPe(t)$ and $\DPp(t)$ (cf.\ e.g.\ \cite[p.~56]{Abe}). 

\begin{proposition} 
\label{HiHe}
The map $\Hi \longrightarrow \He$ given by  
\begin{align*}
b_n \mapsto t_n:=\sum_{|t|=n}t ,
\end{align*} 
and extended as a homomorphism of unital algebras, is an injective co-algebra 
homomorphism. In particular, the Hopf algebra $\Hi$ is a Hopf sub-algebra 
of $\He$. The same is true if $\He$ is replaced by $\Hp$ in these statements.
\end{proposition}

\begin{proof}
The reason why the map is injective is given in the proof of  
Theorem~\ref{HdHat}. 
To prove that $\Hi$ is a sub-co-algebra of $\He$, and of $\Hp$ respectively, 
we need to show that $\DPe(t_n)=\DPp(t_n)=\Di(b_n)$, that is 
\begin{align}
\label{DPe=Di}
\DPe(t_n) = \DPp(t_n) &= \sum_{k=0}^n t_k \otimes t_{n-k}. 
\end{align}
We prove these identities by induction on $n$. 
Since the proof is the same for the two co-products $\DPe$ and $\DPp$, 
we only write it down for $\DPe$. 
For $n=0$ we have $\DPe(\|) = \| \otimes \|$, and for $n=1$ we have 
$\DPe(\Y) = \Y \otimes \| + \| \otimes \Y$. So, formula \eqref{DPe=Di}
holds for $n=0$ and for $n=1$ because $t_0=\|$ and $t_1=\Y$. 
Now we suppose that it holds up to a fixed $n \geq 1$, and we show that 
it holds for $n+1$. 

Using the decomposition $t=r \vee s$ of a tree into its left and right 
components, we have
\begin{align*}
t_{n+1} 
&= \sum_{m=0}^n\underset{|r|=n-m}{\sum_{|s|=m}} r \vee s.  
\end{align*} 
Moreover, by the induction hypothesis, we know that 
\begin{equation*}
\sum_{|s|=m} s\ix1 \otimes s\ix2 = \sum_{|s|=m} \DPe(s) = \DPe(t_m) 
= \sum_{k=0}^m t_{m-k} \otimes t_m. 
\end{equation*}
Using these two identities, we obtain
\begin{align*}
\DPe(t_{n+1}) 
&= \sum_{m=0}^{n} \underset{|r|=n-m}{\sum_{|s|=m}}\DPe(r \vee s) \\ 
&= \| \otimes \sum_{m=0}^{n} \underset{|r|=n-m}{\sum_{|s|=m}} r\vee s 
+\sum_{m=0}^{n} \underset{|r|=n-m}{\sum_{|s|=m}} (r\vee s\ix1) \otimes s\ix2 \\
&= \| \otimes t_{n+1} 
+\sum_{m=0}^{n} \sum_{k=0}^m (t_{n-m}\vee t_{m-k}) \otimes t_k \\ 
&=\sum_{k=0}^{n+1} t_{n-k+1} \otimes t_k. 
\end{align*}
\end{proof}

In \cite{BFqedtree}, it was shown that the co-action $\da$ of $\Ha$ on itself 
can be extended to two co-actions $\dde$ and $\ddp$ of $\Ha$ on $\He$
and $\Hp$,
respectively. These allow one 
to define the semi-direct Hopf algebra 
$\Ha \ltimes \He$ (also called smash Hopf algebra; cf.\ \cite{Majid,Molnar}), 
which represents the renormalization Hopf algebras for 
quantum electrodynamics. 

\begin{corollary}
The semi-direct Hopf algebra $\C(\Gd) \ltimes \Hi$ is a Hopf 
sub-algebra of the QED renormalization Hopf algebra $\Ha \ltimes \He$. 
\end{corollary}

\begin{proof}
This follows from the fact that the maps $a_n \mapsto t_n$ and 
$b_n \mapsto t_n$ induce an injective algebra and co-algebra homomorphism 
from $\Hd \ltimes \Hi$ to $\Hat \ltimes \He$. To prove this,  
we use Theorem~\ref{HdHat}, Proposition~\ref{HiHe}, and we 
extend the co-action 
$\ddif$, which corresponds to $\da$ thanks to Eq.~\ref{da(an)}, to all 
suitable spaces.
\end{proof}


\section{Relation with the renormalization functor}
\label{Rfunctor}

In \cite{BrouderSchmitt} and \cite{VanDerLaan}, it was shown that 
renormalization in quantum 
field theory can be considered as a functor of bi-algebras. 
More precisely, if $\B $ is a bi-algebra and 
$T(\B )^+=\bigoplus_{n\geq 1}T^n(\B )$, 
then the tensor algebra $T(T(\B )^+)$ can be equipped with the structure 
of a bi-algebra which corresponds to the Epstein-Glaser renormalization of 
quantum field theories in the configuration space, cf.~\cite{EpsteinGlaser}. 
It was also shown that the renormalization of scalar quantum fields 
is ruled by the commutative version of this bi-algebra, realized on the 
double symmetric space $S(S(\B )^+)$. 

In this section, we show that $\Hd$ is isomorphic to the Hopf algebra 
obtained from $T(T(\B )^+)$ when $\B $ is the trivial bi-algebra, by
taking a certain quotient.
In the first subsection, we recall the definition of the bi-algebra
$T(T(\B )^+)$ for a generic bi-algebra $\B $, and of its associated 
Hopf algebra. 
Then, in the second subsection, 
we define a bi-algebra $\Bd$, from which one can obtain the Hopf
algebra $\Hd$ as a quotient.
We finally prove the above isomorphism claim in the third subsection.


\subsection{The bi-algebra $T(T(\B )^+)$}

We now introduce a bi-algebra structure on $T(T(\B )^+)$ which occurred
for the first time in \cite{BrouderSchmitt}. However, the reader
should be aware that our presentation here is the opposite of the one in
\cite{BrouderSchmitt}. 
Since we shall have to deal with two distinguished tensor products 
and two co-products, we fix some specific notation to avoid confusions. 

Let $\B $ be a bi-algebra, and let $x$ and $y$ be elements of $\B $.
We denote the product in $\B $ by $x \cdot y$, and we let
$\DB(x)=\sum x\ix1\otimes x\ix2$ denote the co-product in $\B $,
again using Sweedler's notation. 
In $T(\B )^+$, we
replace the tensor symbol by a comma, that is, given
$a=x^1\otimes\cdots \otimes x^n\in T^n(\B )\subseteq T(\B )^+$, we write  
$a=(x^1,\dots,x^n)$ instead. Furthermore, in $T(T(\B )^+)$
we omit tensor symbols, that is,  
given an element $u=a^1\otimes\cdots\otimes a^n
\in T^n(T(\B )^+)\subseteq T(T(\B )^+)$,
with $a^i\in T(\B )^+$, 
we write $u=a^1\dots a^n $ instead.

On $T(T(\B )^+)$ we consider the algebra structure given by the tensor product 
(whose symbol is omitted), and we define a co-product $\Delta$ recursively, 
starting from the co-product $\DB$ on $\B $, as follows. 

First, for any $x\in \B $, we introduce three linear operators 
on $T(T(\B )^+)$, $A_x$, $B_x$, and $C_x$, which correspond to the 
product by $x$ in $\B $, $T(\B )$, and $T(T(\B )^+)$, respectively. 
More precisely, if $a=(x_1,x_2,\dots,x_n) \in T^n(\B ) \subset T^1(T(\B )^+)$, 
we set 
\begin{align*}
A_x(a)&= (x \cdot x_1, x_2,\dots,x_n) \in T^n(\B )\subset T^1(T(\B )^+), \\ 
B_x(a)&= (x, x_1, x_2,\dots,x_n) \in T^{n+1}(\B )\subset T^1(T(\B )^+), \\ 
C_x(a)&= (x)(x_1, x_2,\dots,x_n) \in T^1(\B ) \otimes T^n(\B ) 
\subset T^2(T(\B )^+),
\end{align*}
and if $u=a_1 a_2 \dots a_n \in T^n(T(\B )^+)$, the operators
$A_x$, $B_x$ and $C_x$ act only on $a_1$: 
\begin{align*}
A_x(u)&= A_x(a_1) a_2 \dots a_n \in T^n(T(\B )^+), \\ 
B_x(u)&= B_x(a_1) a_2 \dots a_n \in T^n(T(\B )^+), \\ 
C_x(u)&= C_x(a_1) a_2 \dots a_n \in T^{n+1}(T(\B )^+).
\end{align*}
In particular, any element $a=(x_1,x_2,\dots,x_n)$ of $T(\B )^+$ 
can be written as $a=B_{x_1}(a')$, with $a'=(x_2,\dots,x_n)$. 
Now, the co-product in $T(T(\B )^+)$ is defined recursively 
on the generators as follows: 
\begin{align}
\Delta \big((x)\big) &= \sum (x\ix1)\otimes (x\ix2), \label{DeltaTTB1} \\ 
\Delta (B_x (a)) &= \sum (A_{x\ix1}\otimes B_{x\ix2} + 
B_{x\ix1}\otimes C_{x\ix2}) \Delta a, \label{DeltaTTB}
\end{align}
where $\sum x\ix1\otimes x\ix2= \DB(x)$. 
For example,
\begin{align*}
\Delta\big((x)\big) &= \sum (x\ix1) \otimes (x\ix2),\\
\Delta\big((x,y)\big) &= \sum (x\ix1,y\ix1) \otimes (x\ix2)(y\ix2)
+ \sum (x\ix1\cdot y\ix1) \otimes (x\ix2,y\ix2), \\ 
\Delta\big((x,y,z)\big) &= 
\sum (x\ix1,y\ix1,z\ix1) \otimes (x\ix2)(y\ix2)(z\ix2)
+ \sum (x\ix1,y\ix1\cdot z\ix1)) \otimes (x\ix2)(y\ix2,z\ix2) \\
&+ \sum (x\ix1\cdot y\ix1,z\ix1) \otimes (x\ix2,y\ix2)(z\ix2)
+ \sum (x\ix1\cdot y\ix1\cdot z\ix1) \otimes (x\ix2,y\ix2,z\ix2).
\end{align*}

The co-unit $\varepsilon$ of $T(T(\B )^+)$ is the algebra homomorphism 
$T(T(\B )^+)\to \C$ whose restriction to $T(\B )^+$ is given by
$\varepsilon((x))=\varepsilon_{\B }(x)$, for $x\in \B $, and
$\varepsilon((x_1,\dots,x_n))=0$ for $n\ge1$.

In \cite{BrouderSchmitt}, it was proved that $T(T(\B )^+)$ is 
a bi-algebra, and that one obtains a Hopf algebra structure on 
the quotient $T(T(\B )^+)/\< (x-\varepsilon_\B(x)1) \>$ by the bi-ideal 
generated by $(x-\varepsilon_\B(x)1)$.

\subsection{The bi-algebra $\Bd$}

Let $\Bd = \<a_0,a_1,a_2,\dots\>$ denote the free associative algebra 
on the variables $a_0,a_1,a_2,\dots$. These are the
same variables as those which generate $\Hd$, except that there is
an extra variable $a_0$ which is different from $1$. 
Then $\Bd$ is an associative unital algebra, and the formula of 
Lemma~\ref{Dd(a_n)}, in which the $Q^{(k)}_m(a)$ are the polynomials
from Definition~\ref{defQ} but without the identification of $a_0$
with 1, defines a co-associative co-product $\Dd$ on the 
generators of $\Bd$. 
The first few values of the co-product in $\Bd$ are:
\begin{align*}
\Delta a_0 &= a_0 \otimes a_0, \\ 
\Delta a_1 &= a_0\otimes a_1 + a_1 \otimes a_0^2,\\
\Delta a_2 &= a_0\otimes a_2 + a_1 \otimes (a_0a_1+a_1a_0)
+ a_2\otimes a_0^3,\\
\Delta a_3 &= a_0\otimes a_3 + a_1 \otimes (a_0a_2+a_2a_0+a_1^2)
+ a_2\otimes (a_0^2a_1+a_0 a_1 a_0 + a_1 a_0^2) + a_3\otimes a_0^4.
\end{align*}
Furthermore, we can define a co-unit as the linear map 
$\varepsilon:\Bd \longrightarrow \C$ given by $\varepsilon(1)=1$, 
$\varepsilon(a_0)=1$ and $\varepsilon(a_n)=0$ for $n\ge 1$.
Thus, $\Bd$ becomes an associative bi-algebra. 
However, $\Bd$ is {\em not} a Hopf algebra, 
because the antipode cannot be defined on $a_0$. 

In order to obtain a Hopf algebra, we have to consider the quotient 
$\Bd/\<a_0-\varepsilon(a_0)1\> $ of the bi-algebra $\Bd$ by the bi-ideal 
which identifies $a_0$ with the unit $1$. 
What we obtain is precisely the Hopf algebra $\Hd$ of formal 
diffeomorphisms.  


\subsection{Recursive definition of the co-product $\Dd$}

We are now ready to prove that the Hopf algebra $\Hd$ of formal 
diffeomorphism can be also obtained as the quotient
$T(T(\B )^+)/\< (x-\varepsilon_{\B}(x)1 )\>$ when $\B $ is the trivial 
bi-algebra. 

\begin{theorem}
Let $\B =\C 1$ be the trivial bi-algebra with $1\cdot 1=1$ 
and $\DB(1)=1 \otimes 1$. Denote by $1^{\otimes n}=(1,\dots,1)$ 
the only element (up to a scalar) in $T^n(\B )$. 
Then the algebra homomorphism 
\begin{equation*}
\varphi: T(T(\C 1)^+)/\< (1)-\varepsilon(1)1 \> \longrightarrow \Hd, 
\quad \varphi(1^{\otimes n})=a_{n-1} ,
\end{equation*}
is a Hopf algebra isomorphism. 
\end{theorem}

\begin{proof}
Of course, the relation between $\Bd$ and $\Hd$ is the same as the 
relation between $T(T(\C 1)^+)$ and $T(T(\C 1)^+)/\< (1)-\varepsilon(1)1 \>$, 
and the map $\varphi$ can be lifted up to a map 
$\varphi:T(T(\C 1)^+)\longrightarrow \Bd$. Therefore, it is enough to 
show that $\varphi$ is an isomorphism of associative bi-algebras 
from $T(T(\C)^+)$ to $\Bd$. 

The map $\varphi$ is clearly a bijection, with inverse map 
$\varphi^{-1}(a_n)=1^{\otimes (n+1)}$. It is an algebra homomorphism 
by definition, so it only remains to prove that it is a co-algebra 
homomorphism, that is 
$\Dd(a_n)=(\varphi\otimes\varphi)\Delta (1^{\otimes (n+1)})$. 
Since $\Delta$ on $T(T(\C 1)^+)$ is defined by the recurrence relations 
\eqref{DeltaTTB1} and \eqref{DeltaTTB}, with $x=1$ and $\DB(1)=1\otimes 1$, 
we only need to show that $\Dd$ satisfies the analogous recurrence relations. 
To do this, denote by $A^\varphi$, $B^\varphi$ and $C^\varphi$ the operators 
acting on $\Bd$ which arise canonically from 
the operators $A_1$, $B_1$ and $C_1$ 
on $T(T(\C 1)^+)$ by using the map $\varphi:T(T(\C 1)^+)\longrightarrow \Bd$. 
Explicitly, for any $u\in \Bd$, we have:
\begin{equation}
\label{ABC}
A^\varphi(u)=u, \quad  
B^\varphi(a_n u)=a_{n+1} u, \quad 
C^\varphi(u)= a_0 u.
\end{equation}
We claim that the recurrence relation \eqref{DeltaTTB} 
in $\Bd$ becomes  
\begin{equation}
\label{recursiveDd}
\Dd(a_{n+1}) 
=(A^\varphi\otimes B^\varphi+B^\varphi\otimes C^\varphi) \Dd(a_n).  
\end{equation}
To see this,
we apply $B^\varphi$ to the identity 
$Q^{(k)}_{m-1}(a)=\sum_{l=0}^{m-1} a_l Q^{(k-1)}_{m-1-l}(a)$ 
of Lemma~\ref{recursiveQ}. This gives
\begin{align*}
B^\varphi\big(Q^{(k)}_{m-1}\big) 
&=\sum_{i=0}^{m-1} B^\varphi\big(a_i Q^{(k-1)}_{m-1-i}\big) \\ 
&=\sum_{i=0}^{m-1} a_{i+1} Q^{(k-1)}_{m-1-i}  
=\sum_{i=1}^{m} a_{i} Q^{(k-1)}_{m-i} \\
&=Q^{(k)}_{m} - a_0 Q^{(k-1)}_{m}. 
\end{align*} 
Using this identity for $m-1=n-k$, the value of the operators 
$A^\varphi$, $B^\varphi$ and $C^\varphi$ computed in Eq.~\eqref{ABC}, 
and the Definition~\ref{defQ} of the polynomials $Q^{(k)}_{m}$,
we obtain for the right-hand side of \eqref{recursiveDd}:
\begin{align*}
(A^\varphi\otimes B^\varphi+B^\varphi\otimes C^\varphi) \Dd(a_n)
&=\sum_{k=0}^n a_k \otimes B^\varphi(Q^{(k)}_{n-k})
+\sum_{k=0}^n a_{k+1} \otimes a_0 Q^{(k-1)}_{n-k} \\ 
&=\sum_{k=0}^n a_k \otimes Q^{(k)}_{n-k+1} 
-\sum_{k=0}^n a_k \otimes a_0 Q^{(k-1)}_{n-k+1} 
+\sum_{k=1}^{n+1} a_k \otimes a_0 Q^{(k-1)}_{n-k+1} \\
&=\sum_{k=0}^n a_k \otimes Q^{(k)}_{n-k+1} 
-a_0 \otimes a_0 Q^{(-1)}_{n+1} 
+a_{n+1} \otimes a_0 Q^{(n)}_0 \\ 
&= \sum_{k=0}^{n+1} a_k \otimes Q^{(k)}_{n+1-k} = \Dd(a_{n+1}).  
\end{align*} 
This proves our claim.
\end{proof}

Note that formulae \eqref{ABC} and \eqref{recursiveDd} provide a recursive 
definition for the co-product $\Dd$ on $\Bd$. 


\section{Formal diffeomorphisms in several variables}
\label{variables}

In this section, we show how to extend the construction of the non-commutative 
Hopf algebra of formal diffeomorphisms to series with several variables. 
We content ourselves with considering the case of two variables only, 
as it serves well as an illustration, and as there is no substantial 
difference to the more general several variables case, except that notation 
becomes considerably more cumbersome.

Let $\Gd_2$ denote the set of pairs $(\alpha,\beta)$ of series 
\begin{align*}
\alpha(x,y) &= x 
+ \sum_{n=1}^\infty \sum_{k=1}^\infty \alpha_{nk}\ x^{n+1} y^{k}, \\
\beta(x,y) &= y 
+ \sum_{k=1}^\infty \sum_{n=1}^\infty \beta_{nk}\ x^{n} y^{k+1},
\end{align*}
in two variables $x$ and $y$, with complex coefficients $\alpha_{nk}$ 
and $\beta_{nk}$. 
Consider the composition 
$(\alpha,\beta) \circ (\mu,\nu) := (\alpha(\mu,\nu), \beta(\mu,\nu))$ 
of such pairs, which is defined as the pair of series
\begin{align*}
\alpha(\mu,\nu) &= \mu(x,y)   
+ \sum_{n=1}^\infty \sum_{k=1}^\infty \alpha_{nk}\ 
\mu(x,y)^{n+1} \nu(x,y)^{k}, \\
\beta(\mu,\nu) &= \nu(x,y)  
+ \sum_{k=1}^\infty \sum_{n=1}^\infty \beta_{nk}\ 
\mu(x,y)^{n} \nu(x,y)^{k+1}. 
\end{align*}
The set $\Gd_2$ is a group. Its unit is given by the pair $(x,y)$, while
the inverse pair of a given pair 
can be found with the help the well-known Lagrange--Good formula 
(see e.g.\ \cite{Gessel})
in dimension $2$. 

Therefore, the co-ordinate ring $\C(\Gd_2) $ is a 
commutative Hopf algebra, isomorphic to the polynomial ring 
$\C[a_1,a_2,\dots;b_1,b_2,\dots]$ in two infinite series of variables 
labelled by natural numbers. 

To find a non-commutative version of this Hopf algebra, we proceed as 
for $\Gd$ in Section~\ref{diffeomorphisms}. 
We consider the non-commutative algebra $\C\<a_1,a_2,\dots;b_1,b_2,\dots\>$ 
in two infinite series of variables, and we introduce the generating series
\begin{align*}
A(x,y) &= x 
+ \sum_{n=1}^\infty \sum_{k=0}^\infty a_{nk}\ x^{n+1} y^{k},\\
B(x,y) &= y 
+ \sum_{k=0}^\infty \sum_{n=0}^\infty b_{nk}\ x^{n} y^{k+1}.
\end{align*} 
The co-product is then defined through the double residues 
\begin{align*}
\Delta A(x,y) &= \resz{z_1} \resz{z_2} 
\frac{1}{(z_1-A(x,y))}  \frac{1}{(z_2-B(x,y))} \otimes A(z_1,z_2), \\ 
\Delta B(x,y) &= \resz{z_1} \resz{z_2} 
\frac{1}{(z_1-A(x,y))}  \frac{1}{(z_2-B(x,y))} \otimes B(z_1,z_2). 
\end{align*}
Following the lines of Section~\ref{sec:formdiff}, one can show that
this co-product is co-associative, and, together with the standard
co-unit, it gives $\C\<a_1,a_2,\dots;b_1,b_2,\dots\>$ the structure
of a Hopf algebra which is neither commutative nor co-commutative.

The explicit form of the action of the co-product 
on the generators $a_n$ and $b_n$ 
can be found using the same arguments as those given in 
Section~\ref{diffeo}. We leave the details to the reader.


\end{document}